\documentclass[11pt]{amsart}

\usepackage{amssymb}

\usepackage{amscd,verbatim,
latexsym,amsmath,amsthm}
\usepackage{a4wide}

\usepackage{color}

\numberwithin{equation}{section}
\newtheorem{theorem}{Theorem}[section]
\newtheorem{proposition}[theorem]{Proposition}
\newtheorem{lemma}[theorem]{Lemma}
\newtheorem{corollary}[theorem]{Corollary}
\theoremstyle{definition}
\newtheorem{definition}[theorem]{Definition}

\def\C{{\mathbb C}}
\def\R{{\mathbb R}}
\def\N{{\mathbb N}}
\def\Z{{\mathbb Z}}
\def\bH{\mathbb H}
\def\bs{\backslash}

\def\Re{{\operatorname{Re \,}}}
\def\Im{{\operatorname{Im \,}}}

\def\Gm{\Gamma}
\def\gm{\gamma}
\def\saku{{\bigtriangleup}}

\def\l{\ell}

\def\H{{\mathcal H}}
\def\F{{\mathcal F}}

\def\sfc{{\mathsf c}}

\def\hyp{{\operatorname{hyp}}}
\def\para{{\operatorname{par}}}
\def\diag{{\operatorname{diag}}}

\def\Pr{{\operatorname{Prim}}}

\makeatletter
        \def\eqnarray{% equationarray
                \stepcounter{equation}%
                \let\@currentlabel=\theequation
                \global\@eqnswtrue
                \global\@eqcnt\z@
                \tabskip\@centering
                \let\\=\@eqncr
                $$\halign to \displaywidth\bgroup\@eqnsel\hskip\@centering
                $\displaystyle\tabskip\z@{##}$&\global\@eqcnt\@ne
                \hfil$\displaystyle{{}##{}}$\hfil
                &\global\@eqcnt\tw@$\displaystyle\tabskip\z@{##}$\hfil
                \tabskip\@centering&\llap{##}\tabskip\z@\cr}
\makeatother

\newfont{\Ma}{msam10 scaled\magstep1} %--12pt math symbol
\def\Box{\hfill \mbox{\Ma \symbol{003}}}

\begin{document}

\title[Periods of automorphic forms]
{Dirichlet series constructed from periods of automorphic forms}
\author[Y. Gon]{Yasuro Gon}
\email{ygon@math.kyushu-u.ac.jp}
\address{Faculty of Mathematics\\ Kyushu University\\
Motooka\\ Fukuoka 819-0395\\ Japan}

\date{\today}

\begin{abstract}
We consider certain Dirichlet series of Selberg type, constructed 
from periods of automorphic forms. We study analytic properties of
these Dirichlet series and show that they have analytic continuation to
the whole complex plane. 
\end{abstract}

\thanks{
This work was partially supported by JSPS Grant-in-Aid for Scientific Research (C) no. 23540020
and (C) no. 26400017}

\maketitle

\section{Introduction }

\subsection{Introduction} 

Let $k$ be a fixed natural number, 
$\Gm$ be a co-finite torsion-free discrete subgroup
of $SL(2,\R)$. In this article, we consider certain 
Dirichlet series constructed from periods of automorphic 
forms for $\Gm$ of weight $4k$. 

Let us recall the definition of periods of automorphic forms.
For a hyperbolic element $\gm = 
\Bigl(
\begin{array}{cc}
a & b \\
c & d
\end{array}
\Bigr)
\in SL(2,\R)$, put 
$ Q_{\gm}(z) = cz^2+(d-a)z-b$.  
\begin{definition}[Periods of automorphic forms]
Let $g$ be a weight $4k$ holomorphic 
cusp form for $\Gm$ and $\gm$ be a hyperbolic
element in $\Gm$. The period integral of $g$ over the closed geodesic 
associated to $\gm$ is defined by
\begin{equation}
\alpha_{2k}(\gm,g) = \int_{z_{0}}^{\gm z_{0}} Q_{\gm}(z)^{2k-1} g(z) \, dz.
\end{equation}
This integral does not depend on the choice on the point $z_{0} \in \bH$ and 
the path from $z_{0}$ and $\gm z_{0}$. 
\end{definition}

It is known that these periods 
$\{\alpha_{2k}(\gm,g) \, | \, \gm \in \Gm \mbox{:hyperbolic} \}$ 
determine automorphic form $g$ uniquely, 
and can be expressed by the Petersson
scalar product with relative Poincar\'e series associated to
hyperbolic elements. (Cf. Katok \cite{Ka})
The relative Poincar\'e series have been studied by several 
authors in connection with the problem of construction of 
cusp forms and choosing spanning sets for the space 
of cusp forms $S_{4k}(\Gm)$.

Let $\Pr(\Gm)$ be the set of primitive hyperbolic 
conjugacy classes of $\Gm$. For a hyperbolic element $\gm \in \Gm$, 
put $\ell(\gm)$ be the length of the closed geodesic associated to $\gm$
and $N(\gm) = \exp(\ell(\gm))$.

\begin{definition}[Dirichlet series $\Xi_{\Gm}(s;g)$] \label{def:xi0}
For $g \in S_{4k}(\Gm)$ and $s \in \C$ with $\Re s > 1$,
define
\begin{eqnarray}
 \Xi_{\Gm}(s;g) &:=&
\sum_{\gm \in \Pr (\Gm)} \sum_{m=1}^{\infty} 
\beta_{2k}(\gm,g) \, N(\gm)^{-ms} \\ 
&=& 
\sum_{\gm \in \Pr (\Gm)} 
\beta_{2k}(\gm,g) \, \frac{N(\gm)^{-s}}{1-N(\gm)^{-s}} \nonumber
\label{eq:thx0}
\end{eqnarray}
with
\begin{equation}
\beta_{2k}(\gm,g) = \frac{\overline{\alpha_{2k}(\gm,g)}}{2^{6k-3} \,
\sinh^{2k-1}(\ell(\gm)/2)}.
\end{equation}
This series is absolutely convergent for $\Re s >1$.
\end{definition}
In this article, we investigate analytic properties of $\Xi_{\Gm}(s;g)$. 
Our main result is the following theorem. (Theorem \ref{th:xi})

\begin{theorem} \label{th:x0} Let $\Gm$ be a co-compact torsion-free discrete 
subgroup of $SL(2,\R)$ and $g \in S_{4k}(\Gm)$.  
The function $\Xi_{\Gm}(s;g)$, defined for $\Re s >1$, has the 
analytic continuation as a meromorphic function on the whole 
complex plane. $\Xi_{\Gm}(s;g)$ has at most simple poles located 
at: 
\begin{enumerate}
\item
$s = \frac{1}{2} -j \pm i r_{n} \quad (j \in \{ 0,1 \},  \, n \ge 1)$ when $k=1$, with the residue
\[ \frac{-4(-1)^{j}}{(\pm 2ir_n - j) (\pm 2ir_n - j+1)}
\langle \varphi_{n}^{(1)} \overline{\varphi}_{n}^{(1)}, 
\, g \rangle, \]
\item
$s = \frac{1}{2} -j \pm i r_{n} \quad (j\ge 0, \, n \ge 1)$ when $k \ge 2$, with the residue
\[ 
4 (-1)^{k+j} 
\langle \varphi_{n}^{(k)} \overline{\varphi}_{n}^{(k)}, 
\, g \rangle
\sum_{h=\max(0, \, j-2k+1)}^{j}
\frac{(-1)^h   \binom{2k+h-3}{h} \binom{2k-1}{j-h}}{\prod\limits_{m=0}^{2k-1}
(\pm 2ir_n - j+m)}
.\]
\end{enumerate}
There are no poles other than described as above. 
Here, 
$\{ 1/4+r_n^2 \}_{n=0}^{\infty}$ are eigenvalues of the Laplacian 
$-y^2(\frac{\partial^2}{\partial x^2}+\frac{\partial^2}{\partial y^2})$ acting 
on $L^2(\Gm \backslash \bH)$, and
$\{ \varphi_n \}_{n=0}^{\infty}$ is the orthonormal basis of 
$L^2(\Gm \backslash \bH)$ such that $-y^2(\frac{\partial^2}{\partial x^2}+\frac{\partial^2}{\partial y^2}) 
\varphi_n = (1/4+r_n^2) \varphi_n$. 
Besides, we put $\partial_{2j} = y^{-2j} \frac{\partial}{\partial z} y^{2j}$
and define
\[
\varphi_{n}^{(k)} = 
\Bigl[ \partial_{2k-2}  \cdots \partial_{2} \partial_{0} \Bigr] 
\varphi_{n}, \quad 
\overline{\varphi}_{n}^{(k)} = 
\Bigl[ \partial_{2k-2}  \cdots \partial_{2} \partial_{0} \Bigr] 
\overline{\varphi_{n}}.   
\]
\end{theorem}

To study the Dirichlet series $\Xi_{\Gm}(s;g)$, 
we define certain functions $\Psi_{\Gm}(s;g)$
also by using the periods of automorphic forms. 
We remark that if $k=1$, $\Psi_{\Gm}(s;g)$ is identified with the
first variation of the Selberg zeta function for $\Gm$ in 
the Teichm\"uller space of the Riemann surface 
$\Gm \backslash \bH$. (See \cite{G}) Here, $\bH$ is the upper half plane. 

\begin{definition} \label{def:psi}
For $g \in S_{4k}(\Gm)$ and a fixed 
point $s \in \C$ with $\Re s > 1$,
Put 
\begin{equation}
\Psi_{\Gm}(s;g) := 
\sum_{\gm \in \Pr (\Gm)} 
\frac{\beta_{2k}(\gm,g)}{\ell(\gm)}
\biggl\{ \sum_{j=1}^{2k} p_{j}(s) \, \frac{d}{ds} \log Z_{\gm}^{(j)}(s) \biggr\}, \label{eq:th1}
\end{equation}
and the sum is absolutely convergent.
Here, 
$Z_{\gm}^{(j)}(s)$ :the local higher Selberg zeta function of rank $j$
associated to $\gm$
and the polynomial $p_{j}(s) \in \Z[2s]$
are given by
\begin{eqnarray}
Z_{\gm}^{(j)}(s) 
&:=& \prod_{m=0}^{\infty} \bigl( 1 - N(\gm)^{-(s+m)} \bigr)
^{{j+m-1}\choose{m}}  \quad (j \in \Z) \\
p_j(s)  &:=& (j-1)! {{2k-1}\choose{j-1}}{{2k+j-2}\choose{j-1}}
            \prod_{i=j+1}^{2k}(2s-i). 
\end{eqnarray}  
\end{definition}

The higher Selberg zeta function of rank $j$ ($j \in \Z$), defined by the following 
absolutely convergent Euler product (for $\Re s > 1$)
\begin{equation} 
Z_{\Gm}^{(j)}(s) := \prod_{\gm \in \Pr(\Gm)}Z_{\gm}^{(j)}(s)  = 
\prod_{\gm \in \Pr(\Gm)}
\prod_{m=0}^{\infty} \bigl( 1 - N(\gm)^{-(s+m)} \bigr)
^{{j+m-1}\choose{m}}, 
\end{equation}
is introduced and studied by Kurokawa, 
Wakayama and Hashimoto \cite{KW2, HW}. For $j \in \Z$, this zeta function 
$Z_{\Gm}^{(j)}(s)$ also has a meromorphic continuation to the whole 
complex plane.

We show that analytic properties of $\Xi_{\Gm}(s;g)$ is reduced 
to that of $\Psi_{\Gm}(s;g)$. (Propositions \ref{prop:diff}
and \ref{prop:sum})
Thus we investigate analytic properties of $\Psi_{\Gm}(s;g)$.
We state the result on analytic properties of $\Psi_{\Gm}(s;g)$. 

\begin{theorem} Let $\Gm$ be a co-compact torsion-free discrete 
subgroup of $SL(2,\R)$ and $g \in S_{4k}(\Gm)$. 
The function $\Psi_{\Gm}(s;g)$, defined for $\Re s >1$, has the 
analytic continuation as a meromorphic function on the whole 
complex plane. 
$\Psi_{\Gm}(s;g)$ has at most simple poles located 
at:
 \[ s = \frac{1}{2} \pm i r_{n} \quad (n \ge 1). \nonumber \]
There are no poles other than described as above. 
$\Psi_{\Gm}(s;g)$ satisfies the functional equation
\begin{equation}
\Psi_{\Gm}(1-s;g) = \Psi_{\Gm}(s;g).
\end{equation}
\end{theorem}

This theorem is proved by the resolvent type
trace formulas. (Theorem \ref{th:2})
 
\section{Preliminaries} 

In this section we introduce basic objects and fix notations. 

%\subsection{Notations} 
%We denote by $\N$ the set of natural numbers, i.e.
%$\N=\{1,2,3,\ldots \}$. Put $\N_{0} = \N \cup \{ 0 \}$. 
%The cardinality of a finite set $S$ is denoted by $\# S$. 
 
\subsection{The resolvent of the Laplacian}
For an element 
$\gm = 
\Bigl(
\begin{array}{cc}
a & b \\
c & d
\end{array}
\Bigr)
\in SL(2,\R)$ and
a point $z \in \bH$, put
\[
\gm z = \frac{az+b}{cz+d}, 
\quad  j(\gm,z) = cz+d.
\]
Let $X$ be a Riemann surface of type $(g,n)$ with $2g+n > 2$
and $\Gm$ be a co-finite 
torsion-free discrete subgroup of $SL(2,\R)$ such that
$X \cong \Gm \bs \bH$. Here, 
$\bH=\{z \in \C \, | \, \Im z > 0\}$ is the upper half plane 
with the Poincar\'e metric $y^{-2}(dx^2+dy^2)$. 
The group $\Gm$ is generated by $2g$ 
hyperbolic elements $A_{1}, B_{1},\ldots, A_{g}, B_{g}$
and $n$ parabolic elements $S_{1},\ldots,S_{n}$ 
satisfying the single relation
\[ A_{1}B_{1}A_{1}^{-1}B_{1}^{-1} \cdots A_{g}B_{g}A_{g}^{-1}B_{g}^{-1}
   S_{1} \cdots S_{n} = e. \]

Let $k,\l$ be two integers.  
A smooth complex valued function $f$ on $\bH$ is called
an automorphic form of weight $(2k,2\ell)$ with respect
to the group $\Gm$ if for any $z \in \bH$ and $\gm \in \Gm$,
\[ f(\gm z) = j(\gm,z)^{2k} \overline{j(\gm,z)^{2\l}} 
 f(z). \]
An automorphic form of weight $2k$ is meant for 
that of weight $(2k,0)$.  

We remark that automorphic forms of weight $(2k,2\l)$
correspond to tensors of type $(k,\l)$ on the Riemann surface 
$X \cong \Gm \bs \bH$. We denote by $\H^{k,\l}$ the Hilbert
space of automorphic forms of weight $(2k,2\l)$ with the
scalar product
\begin{equation} \label{eq:inner}
\langle f,g \rangle = 
\int_{\Gm \bs \bH} f(z) \overline{g(z)} y^{2k +2\l}
\, \frac{dxdy}{y^2}.  
\end{equation} 
For each integer $k$ we consider the Laplacian
\begin{equation}
\saku_{k} = \bar{\partial}_{k}^{*} \bar{\partial}_{k}
= - y^{2-2k} \frac{\partial}{\partial z}
y^{2k} \frac{\partial}{\partial \bar{z}}
= - \frac{1}{4} 
\biggl[ y^2 \biggl( \frac{\partial^2}{\partial x^2} 
+ \frac{\partial^2}{\partial y^2} \biggr) 
-2 \sqrt{-1} k y \frac{\partial}{\partial x} \biggr]
\end{equation}
in the Hilbert space $\H^{k}=\H^{k,0}$. 
Here, $\bar{\partial}_{k} = \frac{\partial}{\partial \bar{z}}
= \frac{1}{2} \bigl( \frac{\partial}{\partial x} 
+ \sqrt{-1} \frac{\partial}{\partial y} \bigr)$ is consider 
as an operator from $\H^{k}$ to $\H^{k,1}$, 
and $\bar{\partial}_{k}^{*} = 
-y^{2-2k} \frac{\partial}{\partial z} y^{2k}$
is the adjoint operator to $\bar{\partial}_{k}$
in the scalar product (\ref{eq:inner}), 
acting from $\H^{k,1}$ to $\H^{k}$,
where $\frac{\partial}{\partial z}
= \frac{1}{2} \bigl( \frac{\partial}{\partial x} 
- \sqrt{-1} \frac{\partial}{\partial y} \bigr)$.
The operator $\saku_{k}$ is self-adjoint and non-negative
in $\H^{k}$. We denote by $\Omega^{k}(X) = S_{2k}(\Gm)$ the 
subspace $\ker \saku_{k} = \ker \bar{\partial}_{k}$
in $\H^{k}$, consisting of holomorphic cusp forms of 
weight $2k$. 
%We also put the subspace 
%$\Omega^{k,1}(X) = \ker \bar{\partial}_{k}^{*}
%= \coker \, \bar{\partial}_{k}$ in $\H^{k,1}$, 
%which is the Kodaira-Serre dual of $\Omega^{1-k}(X)$. 

Let us denote by $Q_{s}^{(k)}(z,z')$ the resolvent kernel
of the Laplacian $\saku_{k}$ on the upper half-plane $\bH$,
i.e. $Q_{s}^{(k)}(z,z')$ is the kernel of the operator
$\bigl( \saku_{k}+\frac{1}{4}(s-2k)(s-1) \bigr)^{-1}$ 
for $\Re s \ge 1$.
The kernel $Q_{s}^{(k)}(z,z')$ is smooth for $z \ne z'$
and is holomorphic in $s$ on the whole complex plane. 
The kernel $Q_{s}^{(k)}$ has an important 
property that $Q_{s}^{(k)}(\sigma z, \sigma z') 
= Q_{s}^{(k)}(z,z')$ for any $\sigma \in SL(2,\R)$
and $z,z' \in \bH$. For $k=0$ the kernel 
$Q_{s}^{(k)}$ is given by the explicit formula
\begin{equation} \label{eq:qs0}
Q_{s}^{(0)}(z,z') = \frac{\Gm(s)^2}{\pi \Gm(2s)}
\biggl(
1 - \bigg| \frac{z-z'}{\bar{z}-z'} \biggr|^2 \biggr)^{s}
{}_{2}F_{1} \biggl(s,s;2s;
1- \bigg| \frac{z-z'}{\bar{z}-z'} \biggr|^2 \biggr)
\end{equation}
where ${}_{2}F_{1}(a,b;c;z)$ is the hypergeometric function.

We denote by $G_{s}^{(k)}(z,z')$ the resolvent kernel
of the Laplacian $\saku_{k}$ on the Riemann surface $X$,
i.e. $G_{s}^{(k)}(z,z')$ is the kernel of the operator
$\bigl( \saku_{k}+\frac{1}{4}(s-2k)(s-1) \bigr)^{-1}$ 
on the Riemann
surface $X = \Gm \bs \bH$. For $\Re s > 1$ and $z \ne \gm z',
 \gm \in \Gm$, the kernel $G_{s}^{(0)}$ is given by the
absolute convergent series
\begin{equation} \label{eq:gs0}
G_{s}^{(0)}(z,z') 
= \sum_{\gm \in \Gm} Q_{s}^{(0)}(z,\gm z'),
\end{equation}
which admits term-by-term differentiation with respect to the
variables $z$ and $z'$. 
The kernel $G_{s}^{(0)}(z,z')$ is smooth for 
$z \ne \gm z', \gm \in \Gm$, admits a meromorphic continuation
in $s$ on the whole complex plane.

\section{The function $\Psi_{\Gm}(s;g)$}

\subsection{Poincar\'e series $F_{s}^{(k)}(z)$}

Let $\Gm$ be as in the previous section.
We introduce a certain Poincar\'e series $F_{s}^{(k)}(z)$ constructed from 
a $2k$-th derivative of the resolvent kernel $Q_{s}^{(0)}(z,z')$.

\begin{definition}
Put $\displaystyle{L_{2j} = (yy')^{-2j} \frac{\partial^2}{\partial z \partial z'} (yy')^{2j}}$. 
The function $F_{s}^{(k)}$ on $\bH$ is defined for $\Re s > 1$ by
\begin{equation}
F_{s}^{(k)}(z) := 
L_{2k-2}L_{2k-4} \cdots L_{2}L_{0}
\biggl.   \Bigl( G_{s}^{(0)}(z,z') - Q_{s}^{(0)}(z,z') \Bigr)  
\biggr|_{z'=z}, 
\end{equation}
where $G_{s}^{(0)}$ and $Q_{s}^{(0)}$ are the resolvent kernels 
(\ref{eq:gs0}), (\ref{eq:qs0}) of the
Laplacian on the Riemann surface $X$ and 
on the upper half plane $\bH$ respectively.  
\end{definition}

Since $G_{s}(z,z')$ admits term-by-term differentiation with respect to the
variables $z$ and $z'$, we have
\begin{equation}
F_{s}^{(k)}(z)
= \sum_{\gm \in \Gm \setminus\{ e \}}
L_{2k-2}L_{2k-4} \cdots L_{2}L_{0}
\biggl. \bigl( Q_{s}^{(0)}(z, \gm z') \bigr)
\biggr|_{z'=z}.
\end{equation} 
 From the assumption on $\Gm$, we have 
$\Gm \setminus \{e\} = \Gm_{\hyp} \cup \Gm_{\para}$.
Here, $\Gm_{\hyp}$ and $\Gm_{\para}$ are the set of hyperbolic elements
of $\Gm$ and the set of parabolic elements of $\Gm$ respectively. 
We also define two functions
$H_{s}^{(k)}$ and $P_{s}^{(k)}$ on $H$ for $\Re s>1$ by
\begin{eqnarray}
H_{s}^{(k)}(z)
&:=& \sum_{\gm \in \Gm_{\hyp}}
L_{2k-2}L_{2k-4} \cdots L_{2}L_{0}
\biggl. \bigl( Q_{s}^{(0)}(z, \gm z') \bigr)
\biggr|_{z'=z}, \\
P_{s}^{(k)}(z)
&:=& \sum_{\gm \in \Gm_{\para}}
L_{2k-2}L_{2k-4} \cdots L_{2}L_{0}
\biggl. \bigl( Q_{s}^{(0)}(z, \gm z') \bigr)
\biggr|_{z'=z}.
\end{eqnarray}
By definition, we have $F_{s}^{(k)} = H_{s}^{(k)}+P_{s}^{(k)}$.

We collect fundamental properties of $F_{s}^{(k)}$,
$H_{s}^{(k)}$ and $P_{s}^{(k)}$ 
by using the explicit
formula (\ref{eq:qs0}) for $Q_{s}^{(0)}$.

\begin{proposition} \label{prop:fs}
\begin{enumerate}
\renewcommand{\labelenumi}{(\roman{enumi})}
\item The function $F_{s}^{(k)}(z)$ can be written as
\begin{eqnarray}
F_{s}^{(k)}(z) &=& (-1)^k \frac{1}{\pi} \sum_{\gm \in \Gm \setminus \{e\}} 
          \frac{1}{j(\gm,z)^{2k}} \frac{1}{(z - \gm z)^{2k}}
          \frac{\Gm(s+k)^2}{\Gm(2s)} \nonumber \\
  && \times (r-1)^{2k} r^{s-k} {}_{2}F_{1}(s+k,s+k;2s;r) \label{eq:fs1}
\end{eqnarray}
with $\displaystyle{r=r(z, \gm z) = 
1 - \Bigl|\frac{z - \gm z}{\bar{z} - \gm z}\Bigr|^2}$.
Two functions $H_{s}^{(k)}(z)$ and $P_{s}^{(k)}(z)$ have the same expression by
replacing $\Gm \setminus \{e\}$ with $\Gm_{\hyp}$ or $\Gm_{\para}$.
\item The Poincar\'e series $F_{s}^{(k)}(z)$,  $H_{s}^{(k)}(z)$ 
and $P_{s}^{(k)}(z)$ 
are smooth automorphic forms of weight $4k$
for the Fuchsian group $\Gm$, i.e. 
$F_{s}^{(k)}, H_{s}^{(k)}, P_{s}^{(k)} \in \H^{2k}$. 
\end{enumerate}
\end{proposition}
{\it Proof.} We show that for $F_{s}^{(k)}$. 
The other two case are quite similar. \\ 
(i) \, It is sufficient to show that
\begin{eqnarray}
&& L_{2k-2} \cdots L_{2}L_{0} \, Q_{s}^{(0)}(z,z') \nonumber \\
&& =  \frac{(-1)^k \pi^{-1}}{(z - z')^{2k}}
          \frac{\Gm(s+k)^2}{\Gm(2s)}
          (r-1)^{2k} r^{s-k} {}_{2}F_{1}(s+k,s+k;2s;r) \label{eq:fsk1}
\end{eqnarray}
with $\displaystyle{r=r(z,z') = 1 - \Bigl|\frac{z - z'}{\bar{z} - z'}\Bigr|^2}$.
We prove the formula (\ref{eq:fsk1}) by induction on $k$. 
For a smooth function $f(r) = f(r(z,z'))$ with 
$\displaystyle{r(z, z') = 
1 - \Bigl|\frac{z - z'}{\bar{z} - z'}\Bigr|^2}$, we can easily check that
\begin{eqnarray}
\frac{\partial^2}{\partial z \partial z'} 
\biggl[ \frac{(yy')^{2k}}{(z-z')^{2k}} f(r) \biggr]
&& = -2k \frac{(yy')^{2k}}{(z-z')^{2k+2}}\frac{r+2k}{r}
 - 4k \frac{(yy')^{2k}}{(z-z')^{2k+2}}(1-r) f'(r) 
           \nonumber \\
&& \qquad  - \frac{(yy')^{2k}}{(z-z')^{2k+2}} 
(1-r)^2 \bigl\{ r f''(r) + f'(r) \bigr\}.
\end{eqnarray}

Therefore, (\ref{eq:qs0}): the explicit formula 
of $Q_{s}^{(0)}(r(z,z'))$ gives
\begin{eqnarray}
&& L_{0} \, Q_{s}^{(0)}(r) = 
\frac{\partial^2}{\partial z \partial z'} Q_{s}^{(0)}(r(z,z')) \nonumber \\
&& = - \frac{\pi^{-1}}{(z-z')^2} (r-1)^2
\Bigl\{ r \frac{d^2}{dr^2} + \frac{d}{dr} \Bigr\} 
\sum_{n=0}^{\infty} \frac{\Gm(s+n)^2}{\Gm(2s+n)} \frac{r^{s+n}}{n !} \nonumber \\
&& = - \frac{\pi^{-1}}{(z-z')^2} (r-1)^2 \frac{\Gm(s+1)^2}{\Gm(2s)}
       r^{s-1}  {}_{2}F_{1}(s+1,s+1;2s;r). \label{eq:zqs0}
\end{eqnarray}
So (\ref{eq:fsk1}) of the case $k=1$ is proved. 
Put 
\begin{eqnarray}
f_{k}(r) &:=& (-1)^k \pi^{-1}(1-r)^{2k} r^{s-k}
          \frac{\Gm(s+k)^2}{\Gm(2s)} {}_{2}F_{1}(s+k,s+k;2s;r), \\
D_{k} &:=& -2k \frac{r+2k}{r} -4k(1-r) \frac{d}{dr}
          -(1-r)^2 \Bigl\{ r \frac{d^2}{dr^2} + \frac{d}{dr} \Bigr\} 
\nonumber \\ 
      &=:& D_{k,1}+D_{k,2}+D_{k,3}.
\end{eqnarray}
Then we have only to show that 
\begin{equation}
D_{k} f_{k}(r) = f_{k+1}(r). \label{eq:fsk2}
\end{equation} 
We observe that
\begin{eqnarray}
&& D_{k,1} \, f_{k}(r) \nonumber \\
&&= 2k (-1)^{k+1} \pi^{-1} (1-r)^{2k} (r+2k) r^{s-k-1}
          \frac{\Gm(s+k)^2}{\Gm(2s)} {}_{2}F_{1}(s+k,s+k;2s;r).   
\end{eqnarray}

\begin{eqnarray}
&& D_{k,2} \, f_{k}(r) \nonumber \\
&&= 4k (-1)^{k+1} \pi^{-1}
      \frac{\Gm(s+k)^2}{\Gm(2s)} (1-r) 
\frac{d}{dr} \bigl[ r^{s-k} {}_{2}F_{1}(s-k,s-k;2s;r)  \bigr]
        \nonumber \\
&& = 4k (-1)^{k+1} \pi^{-1} \frac{\Gm(s+k)^2}{\Gm(2s)} (1-r)
    \nonumber \\
&& \quad \times \frac{\Gm(2s)}{\Gm(s-k)^2} \frac{\Gm(s-k)\Gm(s-k+1)}{\Gm(2s)}
    r^{s-k-1} {}_{2}F_{1}(s-k,s-k+1;2s;r) \nonumber \\
&& = 4k (-1)^{k+1} \pi^{-1} (1-r)^{2k} r^{s-k-1} (s-k)
          \frac{\Gm(s+k)^2}{\Gm(2s)} {}_{2}F_{1}(s+k,s+k-1;2s;r).
\end{eqnarray}

\begin{eqnarray}
&& D_{k,3} \, f_{k}(r) \nonumber \\
&&= (-1)^{k+1} \pi^{-1}
      \frac{\Gm(s+k)^2}{\Gm(2s)} (1-r)^2 r^{-1}
\Bigl( r \frac{d}{dr} \Bigr)^2 \bigl[ r^{s-k} {}_{2}F_{1}(s-k,s-k;2s;r)  \bigr]
        \nonumber \\
&&= (-1)^{k+1} \pi^{-1} \frac{\Gm(s+k)^2}{\Gm(2s)} (1-r)^2
    \nonumber \\
&& \quad \times \frac{\Gm(2s)}{\Gm(s-k)^2} \frac{\Gm(s-k+1)^2}{\Gm(2s)}
    r^{s-k-1} {}_{2}F_{1}(s-k+1,s-k+1;2s;r) \nonumber \\
&&= (-1)^{k+1} \pi^{-1} (1-r)^{2k} r^{s-k-1} (s-k)^2
          \frac{\Gm(s+k)^2}{\Gm(2s)} {}_{2}F_{1}(s+k-1,s+k-1;2s;r).
\end{eqnarray}

Therefore, we have
\begin{eqnarray}
&& \biggl\{ (-1)^{k+1} \pi^{-1} (1-r)^{2k} r^{s-k-1} \frac{\Gm(s+k)^2}{\Gm(2s)}
\biggr\}^{-1} D_{k} \, f_{k}(r) \nonumber \\
&& = 2k(r+2k) \, {}_{2}F_{1}(s+k,s+k;2s;r) + 4k(s-k) \, {}_{2}F_{1}(s+k,s+k-1;2s;r)
\nonumber \\ 
&& \quad  +(s-k)^2 \, {}_{2}F_{1}(s+k-1,s+k-1;2s;r).  \label{eq:fsk3}
\end{eqnarray}

We require a lemma for further calculation. 
\begin{lemma}
\begin{eqnarray}
&& 2k(r+2k) \, {}_{2}F_{1}(s+k,s+k;2s;r) + 4k(s-k) \, {}_{2}F_{1}(s+k,s+k-1;2s;r)
\nonumber \\ 
&& \quad  +(s-k)^2 \, {}_{2}F_{1}(s+k-1,s+k-1;2s;r) \nonumber \\
&& = (s+k)^2(1-r)^2 \,  {}_{2}F_{1}(s+k+1,s+k+1;2s;r).
\end{eqnarray}
\end{lemma}
{\it Proof.} To simplify the notation, we write
\begin{equation}
{}_{2}F_{1}(a,b;c;z) \equiv F, \quad 
{}_{2}F_{1}(a \pm 1,b;c;z) \equiv F(a \pm 1), \nonumber
\end{equation}
\begin{equation}
{}_{2}F_{1}(a,b \pm 1;c;z) \equiv F(b \pm 1), \quad 
{}_{2}F_{1}(a,b;c \pm 1;z) \equiv F(c \pm 1). \nonumber
\end{equation}
We use the following contiguous relation (Cf. (9.2.4) in \cite[p.242]{Lebedev}):
\begin{equation} 
(c-a-b)F+a(1-z)F(a+1)-(c-b)F(b-1) = 0. \label{eq:fcont}
\end{equation}

Let $G(r)$ be the left hand side of this lemma to prove. 

By using the formula (\ref{eq:fcont}) for $a=s+k-1,\, b=s+k,\, c=2s$, we have
\begin{eqnarray}
G(r) && =   2k(r+2k) \, F(s+k,s+k;2s;r) 
+ 4k(s-k) \, F(s+k,s+k-1;2s;r)
\nonumber \\ 
&& \qquad  +(s-k)^2 \, F(s+k-1,s+k-1;2s;r) \nonumber \\
&& = 2k(r+2k) \, F(s+k,s+k;2s;r) 
+ 4k(s-k) \, F(s+k,s+k-1;2s;r)
\nonumber \\ 
&& \qquad  +(s-k) \, \Bigl\{ (1-2k) F(s+k-1,s+k;2s;r) \nonumber \\
&& \qquad  +(s+k-1)(1-r) \, F(s+k,s+k;2s;r) \Bigr\}
\nonumber \\
&& = \big\{ 2k(r+2k) +(s+k-1)(s-k)(1-r) \} F(s+k,s+k;2s;r)
\nonumber \\
&& \qquad +(1+2k)(s-k) \, F(s+k,s+k-1;2s;r).
\end{eqnarray}

By using the formula (\ref{eq:fcont}) for $a=s+k,\, b=s+k,\, c=2s$, 
we have
\begin{eqnarray}
G(r) && = \big\{ 2k(r+2k) +(s+k-1)(s-k)(1-r) \} F(s+k,s+k;2s;r)
\nonumber \\
&& \qquad +(1+2k) \Bigl\{ -2k \, F(s+k,s+k;2s;r)
+(s+k)(1-r) \,  F(s+k+1,s+k;2s;r) \Bigr\} \nonumber \\
&& = (s+k)(1-r) \, \Bigl\{ (1+2k) \, F(s+k+1,s+k;2s) \nonumber \\ 
&& \qquad + (s+k-1) \, F(s+k,s+k;2s;r) \Bigr\}. 
\end{eqnarray}

By using the formula (\ref{eq:fcont}) for $a=s+k,\, b=s+k+1,\, c=2s$, 
we have
\begin{eqnarray}
G(r) && = (s+k)(1-r) \Bigl\{ (1+2k) \, F(s+k+1,s+k;2s) \nonumber \\ 
&& \qquad + (s+k-1) \, F(s+k,s+k:2s;r) \Bigr\} \nonumber \\
&& =  (s+k)(1-r) \Bigl\{  (1+2k) \, F(s+k+1,s+k;2s) 
-(2k+1) \, F(s+k,s+k+1;2s;r) \nonumber \\ 
&& \qquad + (s+k)(1-r) \, F(s+k+1,s+k+1;2s;r) \Bigr\} 
\nonumber \\
&& = (s+k)^2 (1-r)^2 \, F(s+k+1,s+k+1;2s;r). 
\end{eqnarray}
It completes the proof of the lemma. $\Box$

Let us complete the proof of Proposition \ref{prop:fs} (i).  
By the above lemma and  (\ref{eq:fsk3}), we have
\begin{eqnarray}
&& \biggl\{ (-1)^{k+1} \pi^{-1} (1-r)^{2k} r^{s-k-1} \frac{\Gm(s+k)^2}{\Gm(2s)}
\biggr\}^{-1} D_{k} \, f_{k}(r) \nonumber \\ 
&& = (s+k)^2 (1-r)^2 \, {}_{2}F_{1}(s+k+1,s+k+1;2s;r).
\end{eqnarray}
At last we have
\begin{eqnarray}
&& D_{k} \, f_{k}(r) \nonumber \\ 
&& =  (-1)^{k+1} \pi^{-1} (1-r)^{2(k+1)} r^{s-k-1}
 \frac{\Gm(s+k+1)^2}{\Gm(2s)}  {}_{2}F_{1}(s+k+1,s+k+1;2s;r) \nonumber \\
&& = f_{k+1}(r).  
\end{eqnarray}
By the assumption of the induction, (\ref{eq:fsk1}) is proved.  
Finally, since $\frac{\partial}{\partial z}(\gm z) = j(\gm,z)^{-2}$, we have
\begin{eqnarray}
F_{s}^{(k)}(z) && = \sum_{\gm \in \Gm \setminus \{e \}} 
L_{2k-2}L_{2k-4} \cdots L_{2}L_{0}
  \biggl.   \Bigl( Q_{s}^{(0)}(z, \gm z') \Bigr)  \biggr|_{z'=z} \nonumber \\
&& = (-1)^k \frac{1}{\pi} \sum_{\gm \in \Gm \setminus \{e \}} 
          \frac{1}{j(\gm,z)^{2k}} \frac{1}{(z - \gm z)^{2k}}
          \frac{\Gm(s+k)^2}{\Gm(2s)} \nonumber \\
&& \quad \times (r-1)^{2k} r^{s-k} {}_{2}F_{1}(s+k,s+k;2s;r)
\end{eqnarray}
with $\displaystyle{r=r(z, \gm z) = 1 
- \Bigl|\frac{z - \gm z}{\bar{z} - \gm z}\Bigr|^2}$ and 
$\displaystyle{L_{2j} = (yy')^{-2j} \frac{\partial^2}{\partial z \partial z'} 
(yy')^{2j}  }$. \\
(ii) \, It is clear from the expression (\ref{eq:fs1}) for $F_{s}^{(k)}$ by
using the following lemma. $\Box$

\begin{lemma} \label{lem:qz}
For $\gm = \Bigl(
\begin{array}{cc}
a & b \\
c & d
\end{array}
\Bigr) \in \Gm$, put $Q_{\gm}(z) = j(\gm,z) \, (z - \gm z)$, i.e. 
\[ Q_{\gm}(z) = cz^2+(d-a)z-b. \] 
Then the polynomial $Q_{\gm}(z)$ satisfies the following formula for 
$\sigma \in \Gm$. 
\[ Q_{\gm}(\sigma z) = j(\sigma,z)^{-2} \, Q_{\sigma^{-1} \gm \sigma}(z). \]
\end{lemma}
{\it Proof.} Note that $j(g_{1}g_{2},z) = j(g_{1}, g_{2}z) j(g_{2},z)$ 
and $(gz-gz')= j(g,z)^{-1} j(g,z')^{-1}(z-z')$ for $g_{1},g_{2},g \in SL(2,\R)$.
By direct calculation, 
\begin{eqnarray}
Q_{\gm}(\sigma z) &=&  j(\gm, \sigma z) (\sigma z - \gm \sigma z) \nonumber \\
&=&  j(\gm \sigma, z) j(\sigma,z)^{-1} \cdot 
     j(\sigma,z)^{-1} j(\sigma, \sigma^{-1} \gm \sigma z)^{-1}
     (z -\sigma^{-1} \gm \sigma z) \nonumber \\
&=&  j(\sigma,z)^{-2} j(\gm \sigma, z)
     \bigl\{ j(\gm \sigma, z) j(\sigma^{-1} \gm \sigma, z) \bigr\}^{-1}
     (z - \sigma^{-1} \gm \sigma z) \nonumber \\
&=&  j(\sigma,z)^{-2} Q_{\sigma^{-1} \gm \sigma}(z).
\end{eqnarray}
So we have the desired formula. $\Box$

\subsection{Inner product formula}

We study the scalar product 
$\langle F_{s}^{(k)}, g \rangle$ in detail.  
We firstly write down this scalar product as the sum over the hyperbolic conjugacy
classes of $\Gm$. Next we investigate the local term appearing in the sum,
more concretely. 

Using an explicit formula of Poincar\'e series $F_{s}^{(k)}$ 
in Proposition \ref{prop:fs}, 
we have the following formula for $\langle F_{s}^{(k)}, g \rangle$.

\begin{lemma} \label{lem:fsg}
Let $\H(\Gm)$ be the set of the $\Gm$-conjugacy classes in $\Gm_{\hyp}$, 
the hyperbolic elements of $\Gm$. We denote $P_{0}$ be the primitive hyperbolic 
element for a given hyperbolic element $P$.    
Put $\F_{0} = \{ z \in \bH \, | \, 1 \le |z| < N(P_{0}) \}$. Here $N(P_{0})$ is the 
norm of $P_{0}$. 
Let $\tau \in SL(2,\R)$ such that 
$\tau P_{0} \tau^{-1} = \diag (N(P_{0})^{1/2}, N(P_{0})^{-1/2})$.
Then  
we have the following formula for $\langle F_{s}^{(k)}, g \rangle$ 
with $g \in S_{4k}(\Gm)$,
\begin{eqnarray}
&& (-1)^k \int_{X} F_{s}^{(k)}(z) \overline{g(z)} y^{4k} \, \frac{dxdy}{y^2} 
\nonumber \\
&& = \sum_{P \in \H(\Gm)} 
\int_{\F_{0}} \frac{N(P)^k}{(N(P)-1)^{2k}}\frac{1}{z^{2k}}
\, f_{s}^{(k)}(r(z,N(P)z)) \, 
\overline{j(\tau^{-1},z)^{-4k}g(\tau^{-1}z)} y^{4k} \, 
\frac{dxdy}{y^2}. \label{eq:fsg1}
\end{eqnarray}
Here, we write 
$\displaystyle{
 F_{s}^{(k)}(z) = 
\sum_{\gm \in \Gm_{\hyp}} \frac{1}{j(\gm,z)^{2k}} \frac{1}{(z-\gm z)^{2k}}
f_{s}^{(k)}(r(z,\gm z))}$ in Proposition \ref{prop:fs}, 
i.e. $f_{s}^{(k)}$ is given by
\[ f_{s}^{(k)}(r(z,\gm z)) = (-1)^k \pi^{-1} \Gm(s+k)^2 \Gm(2s)^{-1}
(r-1)^{2k} r^{s-k}{}_{2}F_{1}(s+k,s+k;2s;r) \] 
with $\displaystyle{r=r(z,\gm z)
= 1 - \Bigl|\frac{z - \gm z}{\bar{z} - \gm z}\Bigr|^2}$. 
\end{lemma}
{\it Proof.} Set $\F = \Gm \bs \bH$ the fundamental domain for $\Gm$ and 
$Z(P)$ be the centralizer of $P$ in $\Gm$. 
\begin{eqnarray}
&& \int_{\F} F_{s}^{(k)}(z) \overline{g(z)} y^{4k} \, \frac{dxdy}{y^2} \nonumber \\
&& = \sum_{P \in \H(\Gm)} \sum_{\sigma \in Z(P) \bs \Gm}
\int_{\F} \frac{1}{j(\sigma^{-1} P \sigma, z)^{2k}} 
\frac{1}{(z-\sigma^{-1}P\sigma z)^{2k}}
\, f_{s}^{(k)}(r(z,\sigma^{-1}P \sigma z)) \overline{g(z)} y^{4k} \, 
\frac{dxdy}{y^2} 
\nonumber \\
&& =  \sum_{P \in \H(\Gm)} \sum_{\sigma \in Z(P) \bs \Gm}
\int_{\F} \frac{1}{j(\sigma, z)^{4k}}\frac{1}{j(P, \sigma z)^{2k}} 
\frac{1}{(\sigma z - P \sigma z)^{2k}}
\, f_{s}^{(k)}(r(\sigma z, P \sigma z)) \nonumber \\
&& \qquad \qquad \qquad \qquad \times
\overline{ j(\sigma, z)^{-4k} g(\sigma z)} y^{4k} \, \frac{dxdy}{y^2} 
\end{eqnarray}
Since $j(\sigma,z)^{-1} \overline{j(\sigma,z)^{-1}} \, y = \Im(\sigma z)$, 
the above expression equals
\begin{eqnarray}
&& \sum_{P \in \H(\Gm)} \sum_{\sigma \in Z(P) \bs \Gm}
\int_{\sigma(\F)} \frac{1}{j(P, z)^{2k}} \frac{1}{(z-Pz)^{2k}}
\, f_{s}^{(k)}(r(z,Pz)) \overline{g(z)} y^{4k} \, \frac{dxdy}{y^2} 
\nonumber \\
&& = \sum_{P \in \H(\Gm)} 
\int_{\F_{P}} \frac{1}{j(P, z)^{2k}} \frac{1}{(z-Pz)^{2k}}
\, f_{s}^{(k)}(r(z,Pz)) \overline{g(z)} y^{4k} \, \frac{dxdy}{y^2}. 
\end{eqnarray}
Here, we set $\F_{P} = \bigcup_{\sigma \in Z(P) \bs \Gm} \sigma(\F)$. 
For a hyperbolic element $P$, there is a $\tau \in SL(2,\R)$ 
such that $\tau P \tau^{-1} = D_{P} := \diag(N(P)^{1/2},N(P)^{-1/2})$
with $N(P)>1$ :the norm of $P$. Then $Pz= \tau^{-1} D_{P} \tau z$ 
and $D_{P}\, \tau z = N(P) \tau z$. $\tau(\F_{P})$ is a fundamental domain of the 
centralizer $Z(D_{P}) = \tau Z(P) \tau^{-1}$. There is a primitive 
element $P_{0}$ such that $P=P_{0}^n$ with $n \in \N$ for $P$. Thus we can replace
$\tau(\F_{P})$ by 
$\F_{0} = \{ z \in \bH \, | \, 1 \le |z| < N(P_{0}) \}$ the fundamental domain
for the group generated by $\tau P_{0} \tau^{-1} = D_{P_{0}}$.
Note that $j(\tau,z)=j(\tau^{-1},\tau z)^{-1}$, then 
the proof is finished. $\Box$

Let us investigate
the local integral corresponding to a hyperbolic conjugacy class 
appearing in the sum (\ref{eq:fsg1}) in Lemma \ref{lem:fsg}.

\begin{definition}
For $P \in \H(\Gm)$ and $g \in S_{4k}(\Gm)$, we define
\begin{equation}
I_{s}^{(k)}(P;g)= \frac{N(P)^k}{(N(P)-1)^{2k}} \int_{\F_{0}} \frac{1}{z^{2k}}
\, f_{s}^{(k)}(r(z,N(P)z)) \overline{g_{\tau}(z)} y^{4k} \, 
\frac{dxdy}{y^2}. 
\end{equation} 
Here, $f_{s}^{(k)}$ is defined in Lemma \ref{lem:fsg} and 
$g_{\tau}(z) = j(\tau^{-1},z)^{-4k}g(\tau^{-1}z)$.
\end{definition}

So we have a paraphrase of Lemma \ref{lem:fsg}, i.e.
\begin{equation}
\langle F_{s}^{(k)}, g \rangle
= \int_{\Gm \bs H} F_{s}(z) \overline{g(z)} y^{4k} \, \frac{dxdy}{y^2} 
= \sum_{P \in \H(\Gm)} I_{s}^{(k)}(P;g).
\end{equation}

We claim that the local term $I_{s}^{(k)}(P;g)$ is the multiple of the 
{\it periods of automorphic forms} over the simple closed geodesic
associated to $P_{0}$.

\begin{proposition} \label{prop:ispg1}
For $P \in \H(\Gm)$, let $P_{0}$ be the primitive element 
such that $P_{0}^{n} = P$ with $n \in \N$.
Then we have
\begin{eqnarray}
I_{s}^{(k)}(P;g) = && - \frac{N(P)^k}{(N(P)-1)^{2k}} 
\frac{N(P_{0})^{k - \frac{1}{2}}}{(N(P_{0})-1)^{2k-1}}
  \, \overline{\alpha_{2k}(P_{0},g)} \nonumber \\
 && \quad \times 
\int_{0}^{\pi}  
 f_{s}^{(k)} \biggl( \frac{4 N(P) \sin^2 \theta}
   {(N(P)-1)^2 \cos^2 \theta + (N(P)+1)^2 \sin^2 \theta} \biggr)
\, \sin^{4k-2} \theta d \theta. \label{eq:ispg11}
\end{eqnarray}
Here, $\alpha_{2k}(P_{0},g) = \int_{z_{0}}^{P_{0}z_{0}}
g(z) Q_{P_{0}}(z)^{2k-1} \, dz$ is the period integral of $g$ over the closed
geodesic associated to $P_{0}$, which does not depend on the choice
of the point $z_{0}$ and the path from $z_{0}$ and $P_{0}z_{0}$. 
The polynomial $Q_{P_{0}}$ is given by 
$Q_{P_{0}}(z) = (z-P_{0}z) \, j(P_{0},z)$. 
\end{proposition}
{\it Proof.} By using the polar coordinate for 
$z=Re^{i \theta} \in \F_{0} = \{ z \in \bH \, | \, 1 \le |z| < N(P_{0}) \}$, 
\begin{eqnarray}
&& I_{s}^{(k)}(P;g)  \nonumber \\
&& \quad = \frac{N(P)^k}{(N(P)-1)^{2k}} \int_{0}^{\pi} \! \! \! \int_{1}^{N(P_{0})}
   R^{-2k}e^{-2 k i \theta}
   f_{s}^{(k)} \biggl( \frac{4 N(P) \sin^2 \theta}
   {(N(P)-1)^2 \cos^2 \theta + (N(P)+1)^2 \sin^2 \theta} \biggr)
   \nonumber \\
&& \qquad \times \overline{g_{\tau}(R e^{i \theta})} 
R^{4k-2} \sin^{4k-2} \theta \, RdR d \theta 
\nonumber \\
&& \quad = \frac{N(P)^k}{(N(P)-1)^{2k}} \int_{0}^{\pi}  
 f_{s}^{(k)} \biggl( \frac{4 N(P) \sin^2 \theta}
   {(N(P)-1)^2 \cos^2 \theta + (N(P)+1)^2 \sin^2 \theta} \biggr)
\nonumber \\
&& \qquad \times \biggl[  \int_{1}^{N(P_{0})}
R^{2k-1} e^{-(2k-1)i \theta} 
\overline{g_{\tau}(R e^{i \theta})} \, e^{-i \theta} dR
\biggr] \, \sin^{4k-2} \theta d \theta.
\end{eqnarray}
Let us consider the complex conjugate of the integral in the square bracket
of the last formula. Put $z_{0} = e^{i \theta}$ and use the fact 
that $(z-N(P_{0})z) \, j(D_{P_{0}},z)= (1-N(P_{0})) N(P_{0})^{-1/2} \, z$, 
then it equals
\[ \int_{z_{0}}^{D_{P_{0}} z_{0}} z^{2k-1} g_{\tau}(z) \, dz  
= \frac{N(P_{0})^{k - \frac{1}{2}}}{(1-N(P_{0}))^{2k-1}}
\int_{z_{0}}^{D_{P_{0}} z_{0}} Q_{D_{P_{0}}}(z)^{2k-1} g_{\tau}(z) 
\, dz. \]  
Recall that $D_{P_{0}} = \diag(N(P_{0})^{1/2}, N(P_{0})^{-1/2})
= \tau P_{0} \tau^{-1}$.  
By using Lemma \ref{lem:qz} on $Q_{P_{0}}$ and fact that $g \in S_{4k}(\Gm)$, 
the differential form satisfies that
\begin{eqnarray}
&& Q_{\tau P_{0} \tau^{-1}}(z)^{2k-1} g_{\tau}(z) \, dz \nonumber \\
&& = Q_{P_{0}}(\tau^{-1}z)^{2k-1} j(\tau^{-1},z)^{4k-2} 
j(\tau^{-1},z)^{-4k} g(\tau^{-1}z) j(\tau^{-1},z)^2 d(\tau^{-1} z) 
\nonumber \\
&& = Q_{P_{0}}(\tau^{-1}z)^{2k-1} g(\tau^{-1}z) d(\tau^{-1}z).
\end{eqnarray} 
So we have
\begin{eqnarray}
&& \int_{z_{0}}^{D_{P_{0}} z_{0}} Q_{\tau P_{0} \tau^{-1}}(z)^{2k-1} 
g_{\tau}(z) \, dz 
 =  \int_{\tau^{-1}z_{0}}^{\tau^{-1} \tau P_{0} \tau^{-1} z_{0}} 
Q_{P_{0}}(z)^{2k-1} g(z) \, dz 
\nonumber \\
&& =  \int_{z_{1}}^{P_{0} z_{1}} Q_{P_{0}}(z)^{2k-1} g(z) \, dz 
\end{eqnarray}
with $z_{1} = \tau^{-1} z_{0}$. 
The rest of the proof follows from the fact that the 
differential form $Q_{P_{0}}(z)^{2k-1} g(z) \, dz$ is holomorphic 
on $H$ and $\langle P_{0} \rangle$-invariant. 
Thus the period $\alpha_{2k}(P_{0},g)$
is well-defined and the proof is finished. $\Box$

Let us calculate the definite integral on $\theta$ in (\ref{eq:ispg11}) 
in Proposition \ref{prop:ispg1}. 
Then we have an explicit formula for $I_{s}^{(k)}(P;g)$ by the following
proposition. 

\begin{proposition} \label{prop:jsp}
Put
\begin{equation}
J_{s}^{(k)}(P) = \int_{0}^{\pi}  
 f_{s}^{(k)} \biggl( \frac{4 N(P) \sin^2 \theta}
   {(N(P)-1)^2 \cos^2 \theta + (N(P)+1)^2 \sin^2 \theta} \biggr)
\, \sin^{4k-2} \theta d \theta. 
\end{equation}
Then we have
\begin{eqnarray}
J_{s}^{(k)}(P)
&& = \frac{(-1)^{k-1}}{2^{4k-2}} \frac{\Gm(2s-1)}{\Gm(2s-2k)}
   (N(P)-1)^{2k-1} N(P)^{-s-k+1} \nonumber \\
&& \times 
{}_{2}F_{1}\biggl(-(2k-1),2k;2-2s; \frac{1}{1-N(P)^{-1}} \biggr). 
\end{eqnarray}
\end{proposition}
{\it Proof.}
By definition,
\[ J_{s}^{(k)}(P) = (-1)^k \frac{1}{\pi} \frac{\Gm(s+k)^2}{\Gm(2s)} \int_{0}^{\pi}
  (r-1)^{2k} r^{s-k} {}_{2}F_{1}(s+k,s+k;2s;r) \, \sin^{4k-2} \theta d \theta
\]
with 
\[ r = r(\theta) =
\frac{4N(P) \sin^2 \theta}{(N(P)-1)^2 \cos^2 \theta +(N(P)+1)^2 \sin^2 \theta}
=\frac{4N(P)}{(N(P)-1)^2 \cot^2 \theta +(N(P)+1)^2} \] 
and 
\[ r-1 = -
\frac{(N(P)-1)^2 (\sin \theta)^{-2}}{(N(P)-1)^2 \cot^2 \theta +(N(P)+1)^2}
. \] 
Substituting $\displaystyle{\cot \theta = \frac{N(P)+1}{N(P)-1} \,t}$, then 
we have
$\displaystyle{ r=r(t)=\frac{4N(P)}{(N(P)+1)^2} \frac{1}{1+t^2}}$,   
\[ -\frac{d \theta}{\sin^2 \theta} =  \frac{N(P)+1}{N(P)-1} \, dt 
\quad \mbox{ and } \quad 
 (r-1)^{2k} \sin^{4k-2} \theta \, d \theta = 
   - \frac{(N(P)-1)^{4k-1}}{(N(P)+1)^{4k-1}} \frac{1}{(1+t^2)^{2k}} \, dt. \] 
Thus we have
\begin{eqnarray}
J_{s}^{(k)}(P) = && (-1)^{k-1} 
\frac{1}{\pi} \frac{(N(P)-1)^{4k-1}}{(N(P)+1)^{4k-1}} \frac{\Gm(s+k)^2}{\Gm(2s)}
\biggl( \frac{4 N(P)}{(N(P)+1)^2} \biggr)^{s-k}     \nonumber \\
&& \times \int_{-\infty}^{\infty} 
\frac{1}{(1+t^2)^{s+k}} \, 
{}_{2}F_{1} \biggl( s+k,s+k;2s; \frac{4 N(P)}{(N(P)+1)^2} \frac{1}{1+t^2}
 \biggr) \, dt. 
\end{eqnarray}
Using the power series expansion of the hypergeometric function, 
we carry out the 
integral term by term: this is permissible by dominated convergence theorem. 

Hence,
\begin{eqnarray}
J_{s}^{(k)}(P) 
= &&  (-1)^{k-1} \frac{1}{\pi} \frac{(N(P)-1)^{4k-1}}{(N(P)+1)^{4k-1}} 
    \biggl( \frac{4 N(P)}{(N(P)+1)^2} \biggr)^{s-k} 
\nonumber \\
&& \times \sum_{n=0}^{\infty} \frac{\Gm(s+k+n)^2}{\Gm(2s+n)}
   \biggl( \frac{4 N(P)}{(N(P)+1)^2} \biggr)^n \frac{1}{n !}
   \int_{-\infty}^{\infty} \frac{dt}{(1+t^2)^{s+n+k}}  
\\
= && (-1)^{k-1} \frac{1}{\pi} \frac{(N(P)-1)^{4k-1}}{(N(P)+1)^{4k-1}} 
    \biggl( \frac{4 N(P)}{(N(P)+1)^2} \biggr)^{s-k} 
\nonumber \\
&& \times \sum_{n=0}^{\infty} \frac{\Gm(s+k+n)^2}{\Gm(2s+n)}
   \biggl( \frac{4 N(P)}{(N(P)+1)^2} \biggr)^n \frac{1}{n !}
   \, \pi ^{\frac{1}{2}} \frac{\Gm(s+n+k-\frac{1}{2})}{\Gm(s+n+k)} 
\\
= && (-1)^{k-1} \pi^{- \frac{1}{2}} \frac{(N(P)-1)^{4k-1}}{(N(P)+1)^{4k-1}} 
    \biggl( \frac{4 N(P)}{(N(P)+1)^2} \biggr)^{s-k} 
\nonumber \\
&& \times \sum_{n=0}^{\infty} \frac{\Gm(s+k-\frac{1}{2}+n) \Gm(s+k+n)}{\Gm(2s+n)}
    \, \frac{1}{n!} \biggl( \frac{4 N(P)}{(N(P)+1)^2} \biggr)^n 
\\
= (-1)^{k-1} && \pi^{- \frac{1}{2}} \frac{(N(P)-1)^{4k-1}}{(N(P)+1)^{4k-1}} 
    \biggl( \frac{4 N(P)}{(N(P)+1)^2} \biggr)^{s-k} \nonumber \\
&& \times \frac{\Gm(s+k-\frac{1}{2}) \Gm(s+k)}{\Gm(2s)}
{}_{2}F_{1} \biggl( s+k-\frac{1}{2},s+k;2s; \frac{4 N(P)}{(N(P)+1)^2} 
 \biggr).
\end{eqnarray}
In the above calculation, we used the well-known formula:
\[ \int_{-\infty}^{\infty} \frac{dt}{(1+t^2)^{\nu+1}}
= 2 \int_{0}^{\frac{\pi}{2}} \cos^{2 \nu } \theta \, d \theta 
= \pi^{\frac{1}{2}} \frac{\Gm(\nu+\frac{1}{2})}{\Gm(\nu+1)}, \quad 
\Re \, \nu > - \frac{1}{2}. \]

We require a lemma on the hypergeometric series:
\begin{lemma}
\begin{eqnarray*}
&& {}_{2}F_{1} \biggl( s+k-\frac{1}{2},s+k;2s; \frac{4 N(P)}{(N(P)+1)^2} 
 \biggr) \nonumber \\
&&  = \biggl( \frac{N(P)+1}{N(P)} \biggr)^{2s+2k-1} 
{}_{2}F_{1}\bigl( 2s+2k-1,2k;2s; N(P)^{-1}\bigr).  
\end{eqnarray*}
\end{lemma}
{\it Proof.}
We use the following formula on quadratic transformations of 
the hypergeometric function (Cf. \cite[(9.6.5), p.251]{Lebedev}): 
\begin{equation} \label{eq:qtf}
{}_{2}F_{1} \Bigl( a,a+\frac{1}{2}; c; z \Bigr)
= \biggl( \frac{1+\sqrt{1-z}}{2} \biggr)^{-2a} \! \!
{}_{2}F_{1} \biggl( 2a, 2a-c+1; c; \frac{1-\sqrt{1-z}}{1+\sqrt{1-z}} \biggr)
\end{equation}
with $|\arg(1-z)| < \pi$. Here, $\sqrt{1-z}$ is meant the branch which 
is positive for real $z$ in the interval $(0,1)$.
Put $a=s+k-\frac{1}{2}$, $c=2s$ and $z =\frac{4 N(P)}{(N(P)+1)^2}$
in (\ref{eq:qtf}), then we have
\begin{eqnarray}
&& {}_{2}F_{1} \biggl( s+k-\frac{1}{2},s+k;2s; \frac{4 N(P)}{(N(P)+1)^2} 
 \biggr) \nonumber \\
&& = \biggl( \frac{1}{2} + \frac{1}{2} \frac{N(P)-1}{N(P)+1} 
\biggr)^{-(2s+k-1)} \! \!
{}_{2}F_{1} \biggl(2s+2k-1,2k;2s; 
\frac{1-(N(P)-1)(N(P)+1)^{-1}}{1+(N(P)-1)(N(P)+1)^{-1}}
\biggr) \nonumber \\
&& = \biggl( \frac{N(P)+1}{N(P)} \biggr)^{2s+2k-1}
{}_{2}F_{1} \Bigl( 2s+2k-1, 2k; 2s; N(P)^{-1} \Bigr). \nonumber 
\end{eqnarray}
This completes the proof. $\Box$

We evaluate the hypergeometric series in the above formula. 
\begin{lemma}
\begin{eqnarray*}
&& {}_{2}F_{1} \bigl( 2s+2k-1,2k;2s; N(P)^{-1} \bigr) \nonumber \\
&&  = (1-N(P)^{-1})^{-2k} \frac{\Gm(2s)}{\Gm(2s+2k-1)} 
\frac{\Gm(2s-1)}{\Gm(2s-2k)} \,
{}_{2}F_{1}\biggl(-(2k-1),2k;2-2s; \frac{1}{1-N(P)^{-1}} \biggr).
\end{eqnarray*}
\end{lemma}
{\it Proof.}
Set $x=N(P)^{-1}$. We have 
\begin{equation}
{}_{2}F_{1} (2s+2k-1,2k;2s;x)  
= (1-x)^{-(4k-1)} {}_{2}F_{1} (-(2k-1),2s-2k;2s;x), \label{eq:dg1}
\end{equation}
by the formula (Cf. \cite[(9.5.3), p.248]{Lebedev}): \\
\begin{equation}
 {}_{2}F_{1} (a,b;c;z) = (1-z)^{c-a-b} {}_{2}F_{1} (c-a,c-b;c;z).
\nonumber 
\end{equation}
Next we use the following formula on linear transformations of 
the hypergeometric function (Cf. \cite[(9.5.8), p.249]{Lebedev}):
\begin{eqnarray}
{}_{2}F_{1} (a,b;c;z) 
&& = (1-z)^{-a} \frac{\Gm(c) \Gm(b-a)}{\Gm(c-a) \Gm(b)}
{}_{2}F_{1} \biggl( a,c-b;1+a-b;\frac{1}{1-z} \biggr) \nonumber \\
&& \quad +(1-z)^{-b} \frac{\Gm(c) \Gm(a-b)}{\Gm(c-b) \Gm(a)}
{}_{2}F_{1} \biggl(c-a,b;1-a+b;\frac{1}{1-z} \biggr). \label{eq:958} 
\end{eqnarray}
Put $a=-(2k-1)$, $b=2s-2k$ and $c=2s$ in (\ref{eq:958}), then we have
\begin{eqnarray}
&& {}_{2}F_{1} (-(2k-1),2s-2k;2s;x) \nonumber \\ 
&& = (1-x)^{2k-1} \frac{\Gm(2s) \Gm(2s-1)}{\Gm(2s+2k-1) \Gm(2s-2k)}
{}_{2}F_{1} \biggl( -(2k-1),2k;2-2s;\frac{1}{1-x} \biggr), \label{eq:dg2}
\end{eqnarray}
By (\ref{eq:dg1}) and (\ref{eq:dg2}), we have the desired formula. 
$\Box$

Let us complete the proof of Proposition \ref{prop:jsp}.

Note that $\Gm(2s+2k-1) = \pi^{-1/2} 2^{2s+2k-2}
\Gm(s+k-\frac{1}{2}) \Gm(s+k) $ by the duplication formula 
of the Gamma function. 
By the above two lemmas, we have 
\begin{eqnarray}
J_{s}^{(k)}(P) 
= &&
(-1)^{k-1} \pi^{- \frac{1}{2}} \frac{(N(P)-1)^{4k-1}}{(N(P)+1)^{4k-1}} 
    \biggl( \frac{4 N(P)}{(N(P)+1)^2} \biggr)^{s-k} 
\nonumber \\
&& \times \frac{\Gm(s+k-\frac{1}{2})\Gm(s+k)}{\Gm(2s)}
\biggl( \frac{N(P)+1}{N(P)} \biggr)^{2s+2k-1} \! \!
\biggl( \frac{N(P)}{N(P)-1} \biggr)^{2k}  
\nonumber \\
&& \times 
\frac{\Gm(2s)}{\Gm(2s+2k-1)} 
\frac{\Gm(2s-1)}{\Gm(2s-2k)} \,
{}_{2}F_{1}\biggl(-(2k-1),2k;2-2s; \frac{1}{1-N(P)^{-1}} \biggr)
\nonumber \\
&& = \frac{(-1)^{k-1}}{2^{4k-2}} \frac{\Gm(2s-1)}{\Gm(2s-2k)}
   (N(P)-1)^{2k-1} N(P)^{-s-k+1} \nonumber \\
&& \times 
{}_{2}F_{1}\biggl(-(2k-1),2k;2-2s; \frac{1}{1-N(P)^{-1}} \biggr). 
\end{eqnarray}

So we get the desired formula. $\Box$

By Proposition \ref{prop:ispg1} and Proposition \ref{prop:jsp},
we obtain an explicit formula for $I_{s}^{(k)}(P;g)$, i.e.
\begin{eqnarray*}
I_{s}^{(k)}(P;g) &=& -\frac{N(P)^k}{(N(P)-1)^{2k}}
             \frac{N(P_{0})^{k-\frac{1}{2}}}{(N(P_{0})-1)^{2k-1}}
             \overline{\alpha_{2k}(P_{0},g)} J_{s}^{(k)}(P) \nonumber \\
&=& \frac{(-1)^k}{2^{4k-2}} \frac{N(P_{0})^{k-\frac{1}{2}}}{(N(P_{0})-1)^{2k-1}}
      \overline{\alpha_{2k}(P_{0},g)}
\frac{\Gm(2s-1)}{\Gm(2s-2k)} \frac{N(P)}{N(P)-1} N(P)^{-s}
\nonumber \\
&& \times   
{}_{2}F_{1}\biggl(-(2k-1),2k;2-2s; \frac{N(P)}{N(P)-1} \biggr).   
\end{eqnarray*}

\begin{theorem} \label{th:a} 
For $g \in S_{4k}(\Gm)$ 
and a fixed point $s \in \C$ with $\Re s >1$,
we have
\begin{equation} 
 \langle F_{s}^{(k)}, g \rangle
 =  \int_{\Gm \bs H} F_{s}^{(k)}(z) \overline{g(z)} 
 \, y^{4k} \frac{dx dy}{y^2} 
=  \sum_{P \in \H(\Gm)} I_{s}^{(k)}(P;g)  \label{eq:tha1}
\end{equation}
and $I_{s}^{(k)}(P;g)$ is given by
\begin{eqnarray}
I_{s}^{(k)}(P;g)
&& = (-1)^k
\frac{\overline{\alpha_{2k}(P_{0},g)}}{2^{6k-3} 
\sinh^{2k-1}(2^{-1}\log N(P_{0}))} \, \frac{\Gm(2s-1)}{\Gm(2s-2k)}
\nonumber \\
&& \quad \times {}_{2}F_{1}\biggl(-(2k-1),2k;2-2s; \frac{N(P)}{N(P)-1} \biggr)
 \, \frac{N(P)}{N(P)-1} N(P)^{-s} \label{eq:tha2} \\
&& = (-1)^k
\frac{\overline{\alpha_{2k}(P_{0},g)}}{2^{6k-3} 
\sinh^{2k-1}(2^{-1}\log N(P_{0}))} \, \frac{\Gm(2s-1)}{\Gm(2s-2k)}
\nonumber \\
&& \quad \times \Bigl[ \sum_{j=0}^{2k-1} 
\frac{ (1-2k)_{j} (2k)_{j} }{(2-2s)_{j} \cdot j!} 
\biggl( \frac{N(P)}{N(P)-1} \biggr)^{j} \Bigr]
 \, \frac{N(P)}{N(P)-1} N(P)^{-s} \label{eq:tha3} 
\end{eqnarray}
The series (\ref{eq:tha1}) is absolutely convergent.
\end{theorem}

We can rewrite the above theorem by using local 
higher Selberg zeta functions of rank $j$ $(1 \le j \le 2k)$.
Then, we find that the function $\Psi_{\Gm}(s;g)$ in 
Definition \ref{def:psi} equals that the inner product 
in the above theorem, i.e.
\begin{equation}
\Psi_{\Gm}(s;g) =  (-1)^k \langle F_{s}^{(k)}, g \rangle. 
\end{equation}

\begin{theorem} \label{th:b}
For $g \in S_{4k}(\Gm)$ 
and a fixed point $s \in \C$ with $\Re s >1$,
we have
\begin{eqnarray}
&& \Psi_{\Gm}(s;g) 
 = (-1)^k \langle F_{s}^{(k)}, g \rangle 
= \sum_{P \in \H(\Gm)} (-1)^k I_{s}^{(k)}(P;g) 
\nonumber \\
&& \qquad = 
\sum_{\gm \in \Pr (\Gm)} 
\frac{\overline{\alpha_{2k}(\gm,g)}}{2^{6k-3} \,\ell(\gm) \, 
\sinh^{2k-1}(\ell(\gm)/2)}
\biggl\{ \sum_{j=1}^{2k} p_{j}(s) \, \frac{d}{ds} \log Z_{\gm}^{(j)}(s) 
\biggr\}.
\label{eq:thb}
\end{eqnarray}
Here, $\ell(\gm) = \log N(\gm)$, 
$Z_{\gm}^{(j)}(s)$, the local higher Selberg zeta function of rank $j$,
and the polynomial $p_{j}(s) \in \Z[2s]$ 
are given by
\begin{eqnarray}
Z_{\gm}^{(j)}(s) 
&=& \prod_{m=0}^{\infty} \big( 1 - N(\gm)^{-(s+m)} \big)
^{{j+m-1}\choose{m}},  \\
p_j(s) &=& (j-1)! {{2k-1}\choose{j-1}}{{2k+j-2}\choose{j-1}}
            \prod_{i=j+1}^{2k}(2s-i).
\end{eqnarray}  
The series (\ref{eq:thb}) is absolutely convergent. 
\end{theorem}
{\it Proof.}  By (\ref{eq:tha1}) and (\ref{eq:tha3}) 
in Theorem \ref{th:a}, we have
\begin{eqnarray}
&& (-1)^k \langle F_{s}^{(k)}, g \rangle \nonumber \\
&& = \sum_{P \in \H(\Gm)} 
\frac{\overline{\alpha_{2k}(P_{0},g)}}{2^{6k-3} \, 
\sinh^{2k-1}(\ell(P_{0})/2)}
\biggl\{ \sum_{j=1}^{2k} p_{j}(s) \, 
\biggl( \frac{N(P)}{N(P)-1} \biggr)^{j} 
\biggr\} N(P)^{-s}. 
\label{eq:thb0}
\end{eqnarray}
Here, we put $\displaystyle{p_j(s) = \frac{\Gm(2s-1)}{\Gm(2s-2k)} 
\frac{(1-2k)_{j-1}(2k)_{j-1}}{(j-1)!(2-2s)_{j-1}}}$.
The relation $\log(1-x) = - \sum_{k=1}^{\infty} 
\frac{x^{k}}{k}$
for $|x|<1$ and the definition of $Z_{P_{0}}(s)$ imply
\begin{equation}
 \log Z_{P_{0}}^{(j)}(s) = 
    - {{j+m-1}\choose{m}}
\sum_{m=0}^{\infty} \sum_{k=1}^{\infty}
       \frac{1}{k} N(P_{0})^{-k(m+s)}. \label{eq:logzp0}
\end{equation}
Next, differentiating (\ref{eq:logzp0}) with respect to $s$, we find that 
\begin{eqnarray}
&& \frac{1}{\log N(P_{0})} \frac{d}{ds} \log Z_{P_{0}}^{(j)}(s)
 = \sum_{m=0}^{\infty} \sum_{k=1}^{\infty} {{j+m-1}\choose{m}} N(P_{0})^{-km -ks}
 \nonumber \\
&& \quad = \sum_{k=1}^{\infty} \frac{1}{(1-N(P_{0})^{-k})^j}
     N(P_{0})^{-ks} 
   =\sum_{k=1}^{\infty} 
   \biggl( \frac{N(P_{0})^{k}}{N(P_{0})^{k}-1} \biggr)^j
     N(P_{0})^{-ks}. \label{eq:thb1}
\end{eqnarray}
We can easily check that 
\begin{eqnarray}
p_j(s) &=& \frac{\Gm(2s-1)}{\Gm(2s-2k)} 
\frac{(1-2k)_{j-1}(2k)_{j-1}}{(j-1)!(2-2s)_{j-1}} \nonumber \\
      &=& (j-1)! {{2k-1}\choose{j-1}}{{2k+j-2}\choose{j-1}}
            \prod_{i=j+1}^{2k}(2s-i). \label{eq:thb2}
\end{eqnarray}
Substituting (\ref{eq:thb1}) and (\ref{eq:thb2}) into  
(\ref{eq:thb0}), we have the desired formula. 
The proof of convergence is assured by the following 
corollary. $\Box$ 

\begin{corollary} \label{th:c}
For $g \in S_{4k}(\Gm)$ and $s \in \C$ with $\Re s >1$, 
we have the following estimate. 
\begin{equation}
 \Bigl|
 \Psi_{\Gm}(s;g) 
\Bigr| 
 \le \frac{1}{2^{4k-2}} \Vert y^{2k} g \Vert_{\infty} \biggl\{
\sum_{j=1}^{2k} |p_{j}(s)| \Bigl. \frac{d}{ds} \log Z_{\Gm}^{(j)}(s) 
\Bigr|_{s = \Re s}
\biggr\}.
\end{equation}
Here, $Z_{\Gm}^{(j)}(s)$ $(j=1,2,\ldots,2k)$ are the higher Selberg zeta 
function of rank $j$ for $\Gm$, 
defined by the following absolutely convergent
Euler products for $\Re s>1$, 
\begin{equation}
Z_{\Gm}^{(j)}(s) = \prod_{\gm \in \Pr(\Gm)} 
\prod_{m=0}^{\infty}(1-N(\gm)^{-(m+s)})^{{{j+m-1}\choose{m}}}.
\end{equation} 
\end{corollary}
{\it Proof.} Firstly we estimate $\alpha_{2k}(\gm,g)$ for $\gm \in \Pr(\Gm)$ and 
$g \in S_{4k}(\Gm)$. Let $\tau \in SL(2,\R)$ such that
$ \tau \gm \tau^{-1} = \diag(N(\gm)^{1/2}, N(\gm)^{-1/2})$. 
Take the point $z_{0} \in \bH$ such that $\tau^{-1}(z_{0}) = i$. 
Since the geodesic connecting $i$ and $i N(\gm)$ is the 
imaginary axis, we have
\begin{eqnarray}
\alpha_{2k}(\gm,g) &=& \int_{z_{0}}^{\gm z_{0}} Q_{\gm}(z)^{2k-1} g(z) \, dz 
= \int_{i}^{iN(\gm)} \bigl\{ 
N(\gm)^{- \frac{1}{2}} - N(\gm)^{\frac{1}2{}} \bigr\}^{2k-1}
  z^{2k-1} g_{\tau}(z) \, dz \nonumber \\
&=&  -2^{2k-1} \sinh^{2k-1} (\ell(\gm)/2)  \int_{1}^{N(\gm)} (iy)^{2k-1} \, 
g_{\tau}(iy) \, d(iy) 
\nonumber \\
&=&  2^{k-1} \sinh^{2k-1} (\ell(\gm)/2)  \int_{1}^{N(\gm)} 
y^{2k} g_{\tau}(iy) \, \frac{dy}{y} 
\end{eqnarray}
Therefore, we have
\begin{eqnarray}
|\alpha_{2k}(\gm,g)| &\le& 2^{2k-1} \sinh^{2k-1}(\ell(\gm)/2) 
\Bigl( \sup_{y>0} |\Im(\tau^{-1}(iy))^{2k} g(\tau^{-1}(iy))| \Bigr)
\int_{1}^{N(\gm)} \frac{dy}{y} \nonumber \\
&\le& 2^{2k-1} \, \ell(\gm) \sinh^{2k-1}(\ell(\gm)/2) 
\Vert y^{2k} g \Vert_{\infty}.
\label{eq:thc1}
\end{eqnarray}
Secondly in (\ref{eq:thb0}), for $s=\sigma+it \in \C$ with $\sigma>1$ and 
a hyperbolic element $P$, 
\begin{equation}
\biggl| p_{j}(s) \biggl( \frac{N(P)}{N(P)-1} \biggr)^{j} N(P)^{-s} \biggr|
= |p_{j}(s)| \biggl( \frac{N(P)}{N(P)-1} \biggr)^{j} 
N(P)^{-\sigma}. \label{eq:thc3} 
\end{equation} 
From (\ref{eq:thc1}) and (\ref{eq:thc3}), 
we complete the proof. $\Box$

%%%%%%%%%%%%%%%%%%%%%%%%%%%%%%%%%%%%%%%%%%%%%%%%%%%%%%%%%%%%%
%%%%%%%%%%%%%%%%%%%%%%%%%%%%%%%%%%%%%%%%%%%%%%%%%%%%%%%%%%%%%
\section{Analytic continuation of $ \Psi_{\Gm}(s;g) $}
%%%%%%%%%%%%%%%%%%%%%%%%%%%%%%%%%%%%%%%%%%%%%%%%%%%%%%%%%%%%%
%%%%%%%%%%%%%%%%%%%%%%%%%%%%%%%%%%%%%%%%%%%%%%%%%%%%%%%%%%%%

Hereafter, we assume that $\Gm$ is {\it co-compact} torsion-free,
i.e. $X$ is a compact Riemann surface of genus $g \ge 2$ . 

\subsection{A variant of the resolvent trace formula}
Since we assume that $\Gm$ is co-compact, the Laplacian $\saku_{0}$ has no
continuous spectrum on $L^{2}(\Gm \bs \bH)$. The eigenvalues of 
\[ 4 \, \saku_{0} = -y^2 \biggl( \frac{\partial^2}{\partial x^2} + \frac{\partial^2}{\partial y^2}   \biggr)
\]
forms a countable subset of non-negative real numbers enumerated as
\[  0 = \lambda_{0} < \lambda_{1} \le \lambda_{2} \le \cdots \le
\lambda_{n} \le \cdots \]
so that each eigenvalues occurs in this sequences with its multiplicity. 

Let $\{ \varphi_{n} \}_{n \ge 0}$ be the orthonormal basis of $L^{2}(\Gm \bs \bH)$
such that $\varphi_{n} \in C^{\infty}(\Gm \bs \bH)$ and 
$4 \, \saku_{0} \varphi_{n} = \lambda_{n} \varphi_{n}$. 
Put $\lambda_{n} =1/4+ r_{n}^2$ for each $n$. 
Recall that $G_{s}^{(0)}$ is the kernel function of the operator
$(\saku_{0} -\frac{1}{4}s(1-s))^{-1}$ on the Riemann surface 
$X = \Gm \bs \bH$.  

\begin{proposition} \label{prop:gs1}
Let $m \in \N$ and $s \in \C$ be such that $m \ge 1$ and $\Re s>1$, 
then the function
$\bigl( - \frac{1}{2s-1}\frac{d}{ds} \bigr)^{m} G_{s}^{(0)}(z,z')$ has 
unique continuous extension to all of $(\Gm \bs \bH)^2$ and 
\begin{equation}
\frac{1}{m!}
\Bigl( - \frac{1}{2s-1}\frac{d}{ds} \Bigr)^{m}
G_{s}^{(0)}(z,z') = 
\sum_{n=0}^{\infty} 
\frac{4}{\{(s-\frac{1}{2})^2+r_{n}^2 \}^{m+1}}
\overline{\varphi_{n}(z)}\varphi_{n}(z'). \label{eq:mgs}
\end{equation}
Here the right hand side of this identity converges uniformly
in $(z,z') \in (\Gm \bs \bH)^2$.
\end{proposition}
{\it Proof.} 
Let $s,a \in \C$, $\Re s >1, \Re a >1$. We have 
(Cf. \cite[Theorem 2.1.2, p.46]{F})
\begin{equation}
 G_{s}^{(0)}(z,z') - G_{a}^{(0)}(z,z') 
 = \sum_{n=0}^{\infty} 
\biggl\{ \frac{4}{(s-\frac{1}{2})^2+r_{n}^2} 
        -\frac{4}{(a-\frac{1}{2})^2+r_{n}^2} \biggr\}
\overline{\varphi_{n}(z)}\varphi_{n}(z').  
\label{eq:s-a}
\end{equation}
By differentiating (\ref{eq:s-a}) with respect to $s$, we have formally
\begin{equation}
\frac{1}{m!}
\Bigl( - \frac{1}{2s-1}\frac{d}{ds} \Bigr)^{m}
G_{s}^{(0)}(z,z') = 
\sum_{n=0}^{\infty} 
\frac{4}{\{(s-\frac{1}{2})^2+r_{n}^2 \}^{m+1}}
\overline{\varphi_{n}(z)}\varphi_{n}(z').
\end{equation}
The series 
\begin{eqnarray}
&& \sum_{n=0}^{\infty} \frac{4}{|(s-\frac{1}{2})^2+r_{n}^2|^{m+1}}
\biggl(  \int_{\Gm \bs \bH} |\varphi_{n}(z)|^{2} \, \frac{dx dy}{y^2} 
\biggr)^{\frac{1}{2} \times 2} \nonumber \\
&& =  \sum_{n=0}^{\infty} \frac{4}{|(s-\frac{1}{2})^2+r_{n}^2|^{m+1}}
\end{eqnarray}
converges absolutely since $\sum_{n \ge 1} \lambda_{n}^{-2}$ converges. 
The proof is finished.
See also Proposition 4.2.1 in \cite{GT}.
$\Box$

\begin{proposition} \label{prop:gs3}
Let $m \in \N$ and $s \in \C$ be such that $m \ge 2k+1$ and $\Re s>1$. 
Put $\partial_{2j} = y^{-2j} \frac{\partial}{\partial z} y^{2j}$
and $\partial'_{2j} = {y'}^{-2j} \frac{\partial}{\partial z'} {y'}^{2j}$
Then
\begin{eqnarray*}
&& \frac{1}{m!}
[\partial_{2k-2} \cdots \partial_{2}\partial_{0}]
[\partial'_{2k-2} \cdots \partial'_{2}\partial'_{0}]
\Bigl( - \frac{1}{2s-1}\frac{d}{ds} \Bigr)^{m}
G_{s}^{(0)}(z,z') \\ 
&& = \sum_{n=0}^{\infty} 
\frac{4}{\{(s-\frac{1}{2})^2+r_{n}^2 \}^{m+1}}
[\partial_{2k-2} \cdots \partial_{2}\partial_{0}] \overline{\varphi_{n}(z)}
\, 
[\partial'_{2k-2} \cdots \partial'_{2}\partial'_{0}] \varphi_{n}(z'). 
\label{eq:zmgs}
\end{eqnarray*}
Here the right hand side of this identity converges uniformly
in $(z,z') \in (\Gm \bs \bH)^2$.
\end{proposition}
{\it Proof.} Operate the differential operators on the 
both sides of (\ref{eq:mgs}), 
we have the desired formula (\ref{eq:zmgs}) formally. 
The series 
\begin{eqnarray}
&& \sum_{n=0}^{\infty} \frac{4}{|(s-\frac{1}{2})^2+r_{n}^2|^{m+1}}
\biggl(  \int_{\Gm \bs \bH} 
\Bigl| [\partial_{2k-2} \cdots \partial_{2}\partial_{0}]
\varphi_{n}(z) \Bigr|^{2} 
\, y^{2k} \frac{dx dy}{y^2} 
\biggr)^{\frac{1}{2} \times 2} \nonumber \\
&& =  \sum_{n=0}^{\infty} \frac{4 \lambda_{n}^{2k}}{|(s-\frac{1}{2})^2+r_{n}^2|^{m+1}}
\end{eqnarray}
converges absolutely. The proof is finished. 
$\Box$

\begin{proposition} \label{prop:qs0}
Let $m \in \N$ and $s \in \C$ be such that $m \ge 2k+1$ and $\Re s>1$. Then
\begin{equation}
\lim_{z' \to z}
\Bigl( - \frac{1}{2s-1}\frac{d}{ds} \Bigr)^{m}
L_{2k-2} L_{2k-4} \cdots L_{2} L_{0} \,
Q_{s}^{(0)}(z,z') = 0, \quad z \in \bH. 
\end{equation}
\end{proposition}
{\it Proof.} Put $q_{s}^{(k)}(z,z') 
= L_{2k-2} L_{2k-4} \cdots L_{2} L_{0} \, Q_{s}^{(0)}(z,z')$. 
By using (\ref{eq:fsk1}) in the proof of Proposition \ref{prop:fs},
\begin{equation}
q_{s}^{(k)}(z,z') = (-1)^k \frac{\pi^{-1}}{(z-z')^{2k}} (r-1)^{2k} \frac{\Gm(s+k)^2}{\Gm(2s)}
       r^{s-k}  {}_{2}F_{1}(s+k,s+k;2s;r)
\end{equation}
with $\displaystyle{ r(z,z') = 1 - \Bigl| \frac{z-z'}{\bar{z}-z'} 
\Bigr|^2 }$ and the formula (it can be derived from
(9.7.5) and (9.7.6) in \cite[p.257]{Lebedev}):
\begin{eqnarray}
&& {}_{2}F_{1}(s+k,s+k;2s;r) \nonumber \\
&& = \frac{\Gm(2s)}{\Gm(s+k)^2}(1-r)^{-2k} 
\sum_{n=0}^{2k-1} \frac{(-1)^n (2k-1-n)! (s-k)_{n}(s-k)_{n}}{n!}(1-r)^n
\nonumber \\
&& \quad - \frac{\Gm(2s)}{\Gm(s-k)^2}
\sum_{n=0}^{\infty} \frac{(s+k)_{n} (s+k)_{n}}{n! (2k+n)!}
\Bigl[ \log(1-r) + 2 \psi(s+k+n) -\psi(n+1) \Bigr. \nonumber \\
&& \quad \quad - \Bigl. \psi(2k+n+1) \Bigr] (1-r)^n
\quad \mbox{ for } |r-1| < 1, \, |\arg(1-r)| < \pi.
\end{eqnarray}
Here, $(\alpha)_{n} = \Gm(\alpha +n)/\Gm(\alpha)$ and 
$\psi(z) = \Gm'(z)/\Gm(z)$. Therefore, we have
\begin{eqnarray}
&& \biggl\{  (-1)^k \frac{\pi^{-1}}{(z-z')^{2k}} r^{s-k} \bigg\}^{-1}
q_{s}^{(k)}(z,z') \nonumber \\
&& = \sum_{n=0}^{2k-1} \frac{(-1)^n (2k-1-n)! (s-k)_{n} (s-k)_{n}}{n!}(1-r)^n
\nonumber \\
&& \quad - \frac{\Gm(s+k)^2}{\Gm(s-k)^2}
\sum_{n=0}^{\infty} \frac{(s+k)_{n} (s+k)_{n}}{n! (2k+n)!}
\Bigl[ \log(1-r) \Bigr. \nonumber \\
&& \quad \quad  \Bigl. +2 \psi(s+k+n) -\psi(n+1) - \psi(2k+n+1) \Bigr] 
(1-r)^{n+2k}. 
\end{eqnarray}
Note that
\begin{eqnarray}
&& r^{s-k} = \{ 1 + (r-1) \}^{s-k} \nonumber \\
&& = \sum_{j=0}^{2k} {{s-k}\choose{j}} (-1)^j (1-r)^j
 + O((1-r)^{2k}) \quad (r \to 1-0),  
\end{eqnarray}
then we have (after some calculation, see also Proposition 3.1.3 in \cite{GT})
\begin{eqnarray}
q_{s}^{(k)}(z,z') &=& (-1)^k \frac{\pi^{-1}}{(z-z')^{2k}}
\Bigl\{ \sum_{j=0}^{2k-1} a_{j}(s)(1-r)^j + b(s)(1-r)^{2k} \log(1-r) \Bigr. 
\nonumber \\ 
&& \Bigl. \qquad  \qquad \qquad + O((1-r)^{2k}) \Bigr\} \qquad (r \to 1-0) 
\end{eqnarray}
with 
\begin{eqnarray}
a_{j}(s) &=& (-1)^j\frac{(2k-1-j)!}{j!} 
  \prod_{i=0}^{j-1} \bigl\{ s^2-s-k^2-(2i-1)k+i(i+1) \bigr\}, \\
b(s) &=& - \frac{1}{(2k)!} \prod_{i=0}^{k-1}
 \bigl\{ s^2-s-i(i+1) \bigr\}^2.
\end{eqnarray}
We find that 
\begin{eqnarray*}
 && a_{j}(1-s) = a_{j}(s), \quad \deg a_{j}(s) = 2j \quad (0 \le j \le 2k-1), \\  
 &&    b(1-s) = b(s), \quad \deg b(s) = 4k. 
\end{eqnarray*} 
So we have 
\begin{equation}
 \Bigl( - \frac{1}{2s-1}\frac{d}{ds} \Bigr)^{j+1}  a_{j}(s) 
=  \Bigl( - \frac{1}{2s-1}\frac{d}{ds} \Bigr)^{2k+1} b(s) = 0.
\end{equation}
Therefore, if $m \ge 2k+1$ then we have
\begin{equation}
\lim_{z' \to z} \Bigl( - \frac{1}{2s-1}\frac{d}{ds} \Bigr)^{m}
q_{s}^{(k)}(z,z') = 0.
\end{equation}
Here we used the fact
\begin{equation}
\lim_{z' \to z} \frac{1}{(z-z')^2}(1-r)^2
= \lim_{z' \to z} \frac{(\bar{z}-\bar{z}')^2}{(\bar{z}-z')^2(z-\bar{z}')^2} = 0. 
\end{equation}
This completes the proof.
$\Box$

\begin{proposition} \label{prop:fs0}
Let $m \in \N$ and $s \in \C$ be such that $m \ge 2k+1$ and $\Re s>1$. Then
\begin{eqnarray}
&& \lim_{z' \to z}
L_{2k-2} L_{2k-4} \cdots L_{2} L_{0} 
\Bigl( - \frac{1}{2s-1}\frac{d}{ds} \Bigr)^{m}
G_{s}^{(0)}(z,z') \nonumber \\ 
&& =  \Bigl( - \frac{1}{2s-1}\frac{d}{ds} \Bigr)^{m}
 \sum_{\gm \in \Gm \setminus \{ e\}} 
L_{2k-2} L_{2k-4} \cdots L_{2} L_{0} 
\biggl. \bigl( Q_{s}^{(0)}(z, \gm z') \bigr)
\biggr|_{z'=z} \nonumber \\
&& = \Bigl( - \frac{1}{2s-1}\frac{d}{ds} \Bigr)^{m}
F_{s}^{(k)}(z).
\end{eqnarray}
\end{proposition}
{\it Proof.} By using Proposition \ref{prop:qs0} and 
interchanging the order of differentiation. 
$\Box$

By Theorem \ref{th:a} and Propositions    
\ref{prop:gs3} and \ref{prop:fs0}, we have the following 
formula. 

\begin{theorem} \label{th:2}
Let $g \in S_{4k}(\Gm)$.
Define $\partial_{2j} = y^{-2j} \frac{\partial}{\partial z} y^{2j}$
and put
\begin{equation}
\varphi_{n}^{(k)} := 
\Bigl[ \partial_{2k-2}  \cdots \partial_{2} \partial_{0} \Bigr] 
\varphi_{n}, \quad 
\overline{\varphi}_{n}^{(k)} := 
\Bigl[ \partial_{2k-2}  \cdots \partial_{2} \partial_{0} \Bigr] 
\overline{\varphi_{n}}.   
\end{equation}
If $m \ge 2k+1$ and $\Re s >1$, then we have
\begin{equation} \label{eq:m-Psi}
 (-1)^k \sum_{n=1}^{\infty} 
\frac{4}{\{(s-\frac{1}{2})^2+r_{n}^2 \}^{m+1}}
\langle \varphi_{n}^{(k)} \overline{\varphi}_{n}^{(k)}, 
\, g \rangle 
= \frac{1}{m !} \Bigl( - \frac{1}{2s-1}\frac{d}{ds} \Bigr)^{m}
\Psi_{\Gm}(s;g)
\end{equation}
with
\begin{equation}
\Psi_{\Gm}(s;g) = 
\sum_{\gm \in \Pr (\Gm)} 
\frac{\overline{\alpha_{2k}(\gm,g)}}{2^{6k-3} \,\ell(\gm) \, \sinh^{2k-1}(\ell(\gm)/2)}
\biggl\{ \sum_{j=1}^{2k} p_{j}(s) \, \frac{d}{ds} \log Z_{\gm}^{(j)}(s) \biggr\}.
\end{equation}
\end{theorem}

\subsection{Analytic continuation of $\Psi_{\Gm}(s;g)$}

We study analytic properties of $\Psi_{\Gm}(s;g)$. By
using Theorem \ref{th:2}, we have the following theorem.

\begin{theorem} \label{th:3}
The function $\Psi_{\Gm}(s;g)$, defined for $\Re s >1$, has the 
analytic continuation as a meromorphic function on the whole 
complex plane. $\Psi_{\Gm}(s;g)$ has at most simple poles located 
at:
\[ s = \frac{1}{2} \pm i r_{n} \quad (n \ge 1). \nonumber \]
There are no poles other than described as above. 
$\Psi_{\Gm}(s;g)$
satisfy the functional equation
\begin{equation}
\Psi_{\Gm}(1-s;g) = \Psi_{\Gm}(s;g).
\end{equation}
\end{theorem}
{\it Proof.}
By using Theorem \ref{th:2}, 
the left hand side of (\ref{eq:m-Psi}) is a meromorphic function of $s \in \C$
and its poles are located at the points $s=1/2 \pm r_n$ with order $m+1$.
Hence, $\Psi_{\Gm}(s;g)$ is a meromorphic function with at simple poles only 
at $s=1/2 \pm r_n$. This completes the proof. $\Box$

%%%%%%%%%%%%%%%%%%%%%%%%%%%%%%%%%%%%%%%%%%%%%%%%%%%%%%%%%
%%%%%%%%%%%%%%%%%%%%%%%%%%%%%%%%%%%%%%%%%%%%%%%%%%%%%%%%%
\section{Dirichlet series $\Xi_{\Gm}(s;g)$}
%%%%%%%%%%%%%%%%%%%%%%%%%%%%%%%%%%%%%%%%%%%%%%%%%%%%%%%%%
%%%%%%%%%%%%%%%%%%%%%%%%%%%%%%%%%%%%%%%%%%%%%%%%%%%%%%%%%

We study analytic properties of the Dirichlet series 
$\Xi_{\Gm}(s;g)$. We show that $\Xi_{\Gm}(s;g)$ are 
related to $\Psi_{\Gm}(s;g)$ and find the relations between them.
As a result, analytic properties of $\Xi_{\Gm}(s;g)$ are derived
from that of $\Psi_{\Gm}(s;g)$. 

\subsection{The difference of $\Psi_{\Gm}(s;g)$}

We consider the difference of $\Psi_{\Gm}(s;g)$. 
 
\begin{proposition} \label{prop:diff}
For $0 \le l \le 2k-1$, $g \in S_{4k}(\Gm)$ and a fixed 
point $s \in \C$ with $\Re s > 1$,
Put 
\begin{eqnarray}
&& \Psi_{\Gm}^{[0]}(s;g) := \Psi_{\Gm}(s;g), \nonumber \\
&& \Psi_{\Gm}^{[l+1]}(s;g) 
:= \frac{1}{2s+l} \bigl\{ \Psi_{\Gm}^{[l]}(s;g) - \Psi_{\Gm}^{[l]}(s+1;g) \bigr\}\quad (0 \le l \le 2k-2). 
\end{eqnarray}
Then, we have
\begin{equation}
\Psi_{\Gm}^{[l]}(s;g) = 
\sum_{\gm \in \Pr (\Gm)} 
\frac{\overline{\alpha_{2k}(\gm,g)}}{2^{6k-3} \,\ell(\gm) \, 
\sinh^{2k-1}(\ell(\gm)/2)}
\biggl\{ \sum_{j=1}^{2k-l} p_{j}^{[l]}(s) \, \frac{d}{ds} 
\log Z_{\gm}^{(j-l)}(s) \biggr\}
\end{equation}
with
\begin{equation}
p_j^{[l]}(s)  = (j-1)! {{2k-1-l}\choose{j-1}}{{2k+j-2-l}\choose{j-1}}
            \prod_{i=j+1}^{2k-l}(2s+l-i).
\end{equation}
\end{proposition}
{\it Proof.} We prove by induction on $l$. 
It is clear for $l=0$. Let
\begin{equation}
F_{\gm}^{[l]}(s) = \sum_{j=1}^{2k-l} p_{j}^{[l]}(s)f_{j-l}(s)
\end{equation} 
with $f_{j}(s) = \frac{d}{ds} 
\log Z_{\gm}^{(j)}(s)$. 
Firstly we note that
\begin{eqnarray*}
&& p_{j}^{[l]}(s) - p_{j}^{[l]}(s+1) \\
&& = (j-1)! {{2k-1-l}\choose{j-1}} {{2k+j-2-l}\choose{j-1}}
     \prod_{i=j+3}^{2k-l-2}(2s+l-i) \\
&& \quad \times \Bigl\{ (2s-2k+2l+1)(2s-2k+2l)  
                        -(2s+l-j+1)(2s+l-j) \Bigr\} \\
&& = (j-1)! {{2k-1-l}\choose{j-1}} {{2k+j-2-l}\choose{j-1}}
     \prod_{i=j+3}^{2k-l-2}(2s+l-i) \\
&& \quad \times \Bigl\{ -2(2k-l-j)(2s) 
   +(2k-l-j)(2k-1+j-3l) \Bigr\}.  
\end{eqnarray*}
Secondly by using the fact: 
\begin{eqnarray}
&& f_{j}(s)-f_{j}(s+1) =  
\frac{d}{ds} \log Z_{\gm}^{(j)}(s)
- \frac{d}{ds} \log Z_{\gm}^{(j-1)}(s) \nonumber \\
&& = \ell(\gm) \sum_{k=1}^{\infty} 
 \biggl( \frac{N(\gm)^k}{N(\gm)^k-1} \biggr)^{j} N(\gm)^{-ks}
     -  \ell(\gm) \sum_{k=1}^{\infty} 
 \biggl( \frac{N(\gm)^k}{N(\gm)^k-1} \biggr)^{j}
     N(\gm)^{-k(s+1)} \nonumber \\
&& = \ell(\gm) \sum_{k=1}^{\infty} 
\biggl( \frac{N(\gm)^k}{N(\gm)^k-1} \biggr)^{j-1} 
     N(\gm)^{-ks} = f_{j-1}(s). 
\end{eqnarray}
Thus we have
\begin{eqnarray}
&& F_{\gm}(s) - F_{\gm}(s+1) \nonumber \\
&& = \sum_{j=1}^{2k-l} \biggl[ 
     p_{j}^{[l]}(s+1) \Bigl\{ f_{j-l}(s) -f_{j-l}(s+1) \Bigr\}
    \biggr. \nonumber \\
&& \quad + (j-1)! (2k-l-j) 
   {{2k-1-l}\choose{j-1}} {{2k+j-2-l}\choose{j-1}}
  \prod_{i=j+3}^{2k-l-2} (2s+l-i) \nonumber \\
&& \biggl. \qquad \times 
\Bigl\{ -2(2s) +(2k-1+j-3l) \Bigr\} f_{j-l}(s) \biggr] \\
&& = \sum_{j=1}^{2k-l} p_{j}^{[l]}(s+1) f_{j-l-1}(s) \nonumber \\
&& \quad + \sum_{j=1}^{2k-l-1} (j-1)! (2k-l-j) 
 {{2k-1-l}\choose{j-1}} {{2k+j-2-l}\choose{j-1}} \nonumber \\
&& \qquad \times \prod_{i=j+3}^{2k-l-2}(2s+l-i) \cdot 
\Bigl\{ -2(2s) + (2k-1+j-3l) \Bigr\} f_{j-l}(s). 
\end{eqnarray}
Let $a_{j}(s)$ be the coefficient function of 
$f_{j-l-1}(s)$ in the last formula. Then we have 
\begin{eqnarray}
&& a_{j}(s) = (j-1)!  
{{2k-1-l}\choose{j-1}} {{2k+j-2-l}\choose{j-1}}
\prod_{i=j+1}^{2k-l}(2s+2+l-i) \nonumber \\
&& \quad +  (j-2)! (2k-l-j+1)  
{{2k-1-l}\choose{j-2}} {{2k+j-3-l}\choose{j-2}} \nonumber \\
&& \qquad \times \prod_{i=j+2}^{2k-l-2}(2s+l-i) \cdot 
\Bigl\{ -2(2s) + (2k-2+j-3l) \Bigr\} \nonumber \\
&& = \prod_{i=j+2}^{2k-l-2}(2s+2+l-i)
 \biggl[ (j-1)!  {{2k-1-l}\choose{j-1}} {{2k+j-2-l}\choose{j-1}}
   (2s+2l-j+1) \biggr. \nonumber \\
&& \quad +  (j-2)! (2k-l-j+1)  
{{2k-1-l}\choose{j-2}} {{2k+j-3-l}\choose{j-2}} \nonumber \\
&& \qquad \biggl. \times \Bigl\{ -2(2s+l) +(2k-2+j-l) \Bigr\} \biggr].
\end{eqnarray}
We show that $a_{j}(s)$ is divisible by $2s+l$.
\begin{eqnarray}
&& a_{j}(s) = \prod_{i=j+2}^{2k-l-2}(2s+2+l-i)
 \biggl[ (2s+l) \biggl\{ 
(j-1)!  {{2k-1-l}\choose{j-1}} {{2k+j-2-l}\choose{j-1}}
\biggr. \biggl. \nonumber \\
&& \qquad  \biggl.  -2 (j-2)! (2k-l-j+1)  
{{2k-1-l}\choose{j-2}} {{2k+j-3-l}\choose{j-2}}
 \biggr\} \nonumber \\
&& \quad +\biggl\{ 
-(j-1) (j-1)!  {{2k-1-l}\choose{j-1}} {{2k+j-2-l}\choose{j-1}}
\biggr. \nonumber \\
&& \qquad +(j-2)! (2k-l-j+1)  
{{2k-1-l}\choose{j-2}} {{2k+j-3-l}\choose{j-2}}
(2k-2+j-l) \biggr\} \biggr] \nonumber \\
&& =  \prod_{i=j+2}^{2k-l-2}(2s+2+l-i)
 \biggl[ (2s+l) \frac{(2k+j-3-l)!}{(j-1)!(2k-j-l-1)!} \biggr]
\nonumber \\
&& = (2s+l) \, (j-1)! 
{{2k-l-2}\choose{j-1}} {{2k+j-l-3}\choose{j-1}}
 \prod_{i=j+1}^{2k-l-1}(2s+l+1-i) \nonumber \\
&& = (2s+l) \, p_{j}^{[l+1]}(s).
\end{eqnarray}
At last we have
\begin{equation}
F_{\gm}(s)-F_{\gm}(s+1) 
= (2s+l) \sum_{j=1}^{2k-l-1} p_{j}^{[l+1]} f_{j-l-1}(s).
\end{equation}
This completes the proof. $\Box$

\begin{lemma} \label{lem:Psi_l}
For $0 \le l \le 2k-1$, we have
\begin{equation}
\Psi_{\Gm}^{[l]}(s;g) = 
\sum_{j=0}^{l} \sfc_j^{[l]}(s) \,
\Psi_{\Gm}(s+j;g)
\end{equation}
with
\begin{equation}
\sfc_j^{[l]}(s) = \frac{(-1)^j \binom{l}{j}}{\prod\limits_{\substack{i=0 \\ i \ne j} }^{l}(2s+j-1+i)}.
\end{equation}
\end{lemma}
{\it Proof.} By the assumption of the induction on $l$,
\begin{eqnarray}
(2s+l) \, \Psi_{\Gm}^{[l+1]}(s;g) &=& \Psi_{\Gm}^{[l]}(s;g) 
- \Psi_{\Gm}^{[l]}(s+1;g) 
\nonumber \\
&=& \sum_{j=0}^{l+1} (-1)^j 
\frac{(2s+j+l)\binom{l}{j} +(2s+j-1)\binom{l}{j-1} }{\prod\limits_{\substack{i=0 \\ i \ne j} }^{l+1}(2s+j-1+i)} 
\Psi_{\Gm}(s+j;g) \nonumber \\
&=& \sum_{j=0}^{l+1} 
\frac{(-1)^j (2s+l)\binom{l+1}{j}}{\prod\limits_{\substack{i=0 \\ i \ne j} }^{l+1}(2s+j-1+i)} 
\Psi_{\Gm}(s+j;g) 
\end{eqnarray}

This completes the proof. $\Box$

By using Theorem \ref{th:3}, Proposition \ref{prop:diff} 
and Lemma \ref{lem:Psi_l}, we have the following theorem.

\begin{theorem} \label{thm:psi-l}
For $0 \le l \le 2k-1$, the function
\begin{equation}
\Psi_{\Gm}^{[l]}(s;g) = 
\sum_{\gm \in \Pr (\Gm)} 
\frac{\overline{\alpha_{2k}(\gm,g)}}{2^{6k-3} \,\ell(\gm) \, 
\sinh^{2k-1}(\ell(\gm)/2)}
\biggl\{ \sum_{j=1}^{2k-l} p_{j}^{[l]}(s) \, \frac{d}{ds} 
\log Z_{\gm}^{(j-l)}(s) \biggr\}
\end{equation}
with
\[
p_j^{[l]}(s)  = (j-1)! {{2k-1-l}\choose{j-1}}{{2k+j-2-l}\choose{j-1}}
            \prod_{i=j+1}^{2k-l}(2s+l-i),
\]
defined for $\Re s >1$, has the 
analytic continuation as a meromorphic function on the whole 
complex plane. $\Psi_{\Gm}^{[l]}(s;g)$ has at most simple poles located 
at: 
\[ s = \frac{1}{2} - j \pm i r_{n} \quad (0 \le j \le l, \, n \ge 1), \]
and its residue at $s = 1/2 - j \pm i r_{n}$ is given by
\begin{eqnarray}
&& \frac{(-1)^j \binom{l}{j}}{\prod\limits_{\substack{m=0 \\ m \ne j} }^{l}
(\pm 2ir_n - j+m)} \mathrm{Res}_{s=1/2\pm i r_n} \Psi_{\Gm}(s;g)  \nonumber \\
&& = 4 (-1)^k \frac{(-1)^j \binom{l}{j}}{\prod\limits_{m=0}^{l}
(\pm 2ir_n - j+m)}
\langle \varphi_{n}^{(k)} \overline{\varphi}_{n}^{(k)}, 
\, g \rangle 
\end{eqnarray}
There are no poles other than described as above. 
$\Psi_{\Gm}^{[l]}(s;g)$
satisfy the functional equation
\begin{equation}
\Psi_{\Gm}^{[l]}(1-l-s;g) = \Psi_{\Gm}^{[l]}(s;g).
\end{equation}
\end{theorem}

Next we introduce certain functions $\Psi_{\Gm}^{[2k-1,p]}(s;g) $
for $1 \le p \le 2k-1$.

\begin{proposition} \label{prop:sum}
For $g \in S_{4k}(\Gm)$ and a fixed 
point $s \in \C$ with $\Re s > 1$,
Put 
\begin{eqnarray}
 \Psi_{\Gm}^{[2k-1,0]}(s;g) &:=& \Psi_{\Gm}^{[2k-1]}(s;g), \nonumber \\
 \Psi_{\Gm}^{[2k-1,p]}(s;g) 
&:=& \sum_{j=0}^{\infty} \Psi_{\Gm}^{[2k-1,p-1]}(s+j;g)  \nonumber \\
& =& \sum_{j=0}^{\infty} \binom{p+j-1}{j} \Psi_{\Gm}^{[2k-1]}(s+j;g)
\quad (p \ge 1).  \label{eq:Psi-sum-j}
\end{eqnarray}
Then, we have
\begin{equation} \label{eq:Psi-sum}
\Psi_{\Gm}^{[2k-1,p]}(s;g) = 
\sum_{\gm \in \Pr (\Gm)} 
\frac{\overline{\alpha_{2k}(\gm,g)}}{2^{6k-3} \,\ell(\gm) \, 
\sinh^{2k-1}(\ell(\gm)/2)} \frac{d}{ds} 
\log Z_{\gm}^{(2-2k+p)}(s).
\end{equation}
Besides, 
$\Psi_{\Gm}^{[2k-1,p]}(s;g)$ $(1 \le p \le 2k-1)$ has the analytic continuation 
as a meromorphic function on the whole complex plane
and has at most simple poles located at: 
\[
s = \frac{1}{2} -j \pm i r_{n} \quad (j \in \{ 0\}\cup \N, \, n \in \N),
\]
and its residue at $s = 1/2 - j \pm i r_{n}$ is given by
\begin{equation}
4 (-1)^{k+j} 
\langle \varphi_{n}^{(k)} \overline{\varphi}_{n}^{(k)}, 
\, g \rangle
\sum_{h=\max(0, \, j-2k+1)}^{j}
\frac{(-1)^h   \binom{p+h-1}{h} \binom{2k-1}{j-h}}{\prod\limits_{m=0}^{2k-1}
(\pm 2ir_n - j+m)}.
\end{equation}
There are no poles other than described as above. 
\end{proposition}
{\it Proof.} For a primitive hyperbolic $\gm \in \Gm$, 
we have 
\begin{eqnarray*}
&& \frac{1}{\ell(\gm)} \sum_{j=0}^{\infty} 
\frac{d}{ds} \log Z_{\gm}^{(m)}(s+j) 
= \sum_{j=0}^{\infty} \sum_{k=1}^{\infty}  
\biggl( \frac{N(\gm)^{k}}{N(\gm)^{k}-1} \biggr)^m
    N(\gm)^{-k(s+j)} \\
&& = \sum_{k=1}^{\infty}  
\biggl( \frac{N(\gm)^{k}}{N(\gm)^{k}-1} \biggr)^{m+1}
     N(\gm)^{-ks} 
=  \frac{1}{\ell(\gm)} \frac{d}{ds} \log Z_{\gm}^{(m+1)}(s).
\end{eqnarray*}
Therefore, we have (\ref{eq:Psi-sum}). 
We can check that $\Psi_{\Gm}^{[2k-1,p]}(s;g)$ is absolutely convergent for $\Re s>1$
by the expression (\ref{eq:Psi-sum}). 
By Theorem \ref{thm:psi-l}, $\Psi_{\Gm}^{[2k-1,p]}(s;g)$ has the analytic continuation 
as a meromorphic function on the whole complex plane.
$\Psi_{\Gm}^{[2k-1,p]}(s;g)$ has at most simple poles at 
$\displaystyle{s = \frac{1}{2} -j \pm i r_{n}}$ $(j\ge 0, n \ge 1)$
by (\ref{eq:Psi-sum-j}).
We completes the proof. $\Box$

\subsection{Dirichlet series $\Xi_{\Gm}(s;g)$}
Finally we have a paraphrase of analytic properties of 
$\Psi_{\Gm}(s;g)$ in terms of 
the Dirichlet series $\Xi_{\Gm}(s;g)$.

\begin{definition} \label{def:xi}
For $g \in S_{4k}(\Gm)$ and $s \in \C$ with $\Re s > 1$,
define
\begin{eqnarray}
 \Xi_{\Gm}(s;g) &=&
\sum_{\gm \in \Pr (\Gm)} \sum_{m=1}^{\infty} 
\beta_{2k}(\gm,g) \, N(\gm)^{-ms} \nonumber \\ 
&=& 
\sum_{\gm \in \Pr (\Gm)} 
\beta_{2k}(\gm,g) \, \frac{N(\gm)^{-s}}{1-N(\gm)^{-s}}
\end{eqnarray}
with
\[
\beta_{2k}(\gm,g) = \frac{\overline{\alpha_{2k}(\gm,g)}}{2^{6k-3} \,
\sinh^{2k-1}(\ell(\gm)/2)}.
\]
\end{definition}

\begin{theorem} \label{th:xi}
The function $\Xi_{\Gm}(s;g)$, defined for $\Re s >1$, has the 
analytic continuation as a meromorphic function on the whole 
complex plane. $\Xi_{\Gm}(s;g)$ has at most simple poles located at: 
\begin{enumerate}
\item
$s = \frac{1}{2} -j \pm i r_{n} \quad (j \in \{ 0,1 \},  \, n \ge 1)$ when $k=1$, with the residue
\[ \frac{-4(-1)^{j}}{(\pm 2ir_n - j) (\pm 2ir_n - j+1)}
\langle \varphi_{n}^{(1)} \overline{\varphi}_{n}^{(1)}, 
\, g \rangle, \]
\item
$s = \frac{1}{2} -j \pm i r_{n} \quad (j\ge 0, \, n \ge 1)$ when $k \ge 2$, with the residue
\[ 
4 (-1)^{k+j} 
\langle \varphi_{n}^{(k)} \overline{\varphi}_{n}^{(k)}, 
\, g \rangle
\sum_{h=\max(0, \, j-2k+1)}^{j}
\frac{(-1)^h   \binom{2k+h-3}{h} \binom{2k-1}{j-h}}{\prod\limits_{m=0}^{2k-1}
(\pm 2ir_n - j+m)}
.\]
\end{enumerate}
There are no poles other than described as above. 
\end{theorem}
{\it Proof.} Put $p=2k-2$ in Proposition \ref{prop:sum} and 
note that 
\[ \frac{1}{\ell(\gm)} 
\frac{d}{ds} \log Z_{\gm}^{(0)}(s)
 = \sum_{k=1}^{\infty} N(\gm)^{-ks} = \frac{N(\gm)^{-s}}{1-N(\gm)^{-s}}. \]
Then we have
\[ \Xi_{\Gm}(s;g) = \Psi_{\Gm}^{[2k-1,2k-2]}(s;g). \]
The proof is finished. $\Box$

\end{document}